\documentclass[12pt]{amsart}
\usepackage{amsmath,amssymb,amsthm,enumerate}
\usepackage[]{latexsym}

\newcommand{\Z}{{\mathbb Z}}
\newcommand{\Q}{{\mathbb Q}}
\newcommand{\C}{{\mathbb C}}
\newcommand{\CC}{{\mathcal C}}
\newcommand{\R}{{\mathbb R}}

\renewcommand{\P}{{\mathbb P}}
\newcommand{\p}{{\mathfrak p}}
\newcommand{\HH}{{\mathbb H}}
\newcommand{\ds}{\displaystyle}

\newcommand{\D}{{\mathcal D}}

\newcommand{\I}{{\mathcal I}}
\newcommand{\M}{{\mathcal M}}
\newcommand{\FF}{{\mathcal F}}
\newcommand{\cP}{{\mathcal P}}

\newcommand{\Spec}{\operatorname{Spec}}

\newcommand{\q}{{\theta}}

\newcommand{\w}{{\omega }}

\newtheorem{theorem}{Theorem}[section]
\newtheorem{definition}[theorem]{Definition}

\newtheorem{proposition}[theorem]{Proposition}

\newtheorem{corollary}[theorem]{Corollary}

\newtheorem{lemma}[theorem]{Lemma}

\newtheorem{remarks}[theorem]{Remarks}
\newtheorem{example}[theorem]{Example}

\include{thebibliography}

\begin{document}

\title[Jacobian Nullwerte, Periods and Symmetric Equations]{Jacobian Nullwerte, Periods and Symmetric Equations for Hyperelliptic Curves}
\author{ Jordi Gu\`ardia }
\thanks{Partially supported by MCYT BFM2003-06768-C02-02.}
\address{J. Gu\`{a}rdia. Departament de Matem\`{a}tica Aplicada IV.
Es\-co\-la Po\-li\-t\`{e}c\-nica Superior d'Enginyeria dede Vilanova i la Geltr\'{u},
Avinguda V\'{\i}ctor Balaguer s/n. E-08800 Vilanova i la Geltr\'{u}}%
\email{guardia@ma4.upc.es}%

\date{\today}

\begin{abstract}
We propose a solution to the hyperelliptic Schottky  problem, based on the use of Jacobian Nullwerte and symmetric models for hyperelliptic curves. Both ingredients are interesting on its own, since the first provide period matrices which can be  geometrically described, and the second have remarkable arithmetic properties.  \\


\end{abstract}

\maketitle

\tableofcontents

\section*{Introduction}

The problem of determining a complex abelian variety from its period lattice is very-well understood from the theoretical viewpoint. The situation is slightly different when the problem is considered with  a computational insight. Efficient numerical algorithms to find equations of an elliptic curve from its period lattice were available many years ago, and there
 are large tables of elliptic curves both in printed and electronic form (\cite{Antwerp}, \cite{Cremona}, \cite{Magma}). Concerning higher-dimensional abelian varieties, last years have seen a significant progress which has led to the elaboration of tables of hyperelliptic curves whose jacobian variety has a given period lattice (\cite{Wamelen}, \cite{Weng}, \cite{Guardia}, \cite{GGG}). There are two main directions in the proposed solutions.

The first solution is essentially due to Mestre \cite{Mestre}, who considered the case of abelian surfaces. Weber \cite{Weber} generalized his work to Jacobian varieties of hyperelliptic curves of any genus.
The general outline of Mestre's method is the following: given a normalized period matrix in the Siegel upper half space corresponding to an abelian variety, one
 calculates, by means of {\it Thetanullwerte},   certain algebraic invariants of a curve whose Jacobian variety is isomorphic to the  desired abelian variety.  One  deduces from these invariants the field of moduli of the curve and finds an equation of the curve over its field of definition. The computations require a certain degree of accuracy,  tend to produce huge intermediate results,  and yield final equations with large coefficients, which must be reduced by some additional method. Apart from computational issues, this method has a second drawback: its geometric nature overpasses the arithmetic of the problem: the initial abelian variety and the Jacobian of the found curve may be only isomorphic over the algebraic closure of its field of definition.

We proposed a second solution for genus two curves inspired in the use of Jacobian Nullwerte. The ideas in \cite{Guardia} drove us to an algorithm which, given the period lattice of a basis of algebraic differential forms of an abelian surface, finds the equation of a genus two hyperelliptic curve defined over the same field as the differential forms.
With our method the arithmetic is preserved, but it requires a better knowledge of the abelian variety. It has been applied satisfactorily  to build a table of 2-dimensional factors of certain modular Jacobian varieties \cite{GGG} and provide examples of abelian surfaces with several polarizations \cite{GGR}. Unfortunately, the algorithm cannot be applied when we only know a normalized period matrix in the Siegel upper half plane,
since in general these periods do not correspond to algebraic differential forms.

Our initial motivation for the present work was to overcome this difficulty controlling the arithmetic of the problem.
Classical ideas already found in \cite{Weil} explain a method to determine a basis of algebraic differential forms from a normalized period matrix. Jacobian Nullwerte are the key tool for making this construction explicit. The combination of this ideas with our original algorithm for abelian surfaces led us to a particular kind of equation for hyperelliptic curves, which we have called {\em symmetric models}. We studied symmetric equations for elliptic curves and its applications in class field theory in \cite{GuToVe}. We present here the geometric study of symmetric models for hyperelliptic curves of any genus. We describe their arithmetic properties, as well as its interest in relation with the problem mentioned above.

After a very short summary  of basic facts on hyperelliptic curves and their Jacobians, we develop the study of symmetric equations in section \ref{seccio-simetric}.
The next section is devoted to recall the main results concerning  Jacobian Nullwerte for hyperelliptic curves. In section \ref{Thomae-seccio} we recall classical formulas of Thomae relating Thetanullwerte and Jacobian Nullwerte to  Weierstrass points of hyperelliptic curves. Some remarks on the theoretical implications of these results are collected in sections \ref{seccio-remarks-igusa} and \ref{seccio-jacobi-remarks}.
We explain the construction of algebraic differential forms from normalized period matrices in section   \ref{matrius-algebraiques}. We then give in section \ref{Nullwerte-Weierstrass} a general method to find a symmetric equation for a general hyperelliptic curve given a normalized period matrix for it. In the last two sections we particularize the results of the paper for hyperelliptic curves of genus 2 and 3, in which some improvements can be obtained.

\section{Preliminaries on hyperelliptic curves and their Jacobians}

We introduce here the notation which will be used along the paper.
Consider a hyperelliptic curve in Weierstrass form:
$$
C: Y^2=f(X)=(X-\alpha_1)\cdots(X-\alpha_{2g+2}),
$$
so that $W_1=(\alpha_1,0),\dots, W_{2g+2}=(\alpha_{2g+2},0)$ are its
Weierstrass points. We denote  by
$\{\w_j=\frac{x^jdx}y\}_{j=0,\dots,g-1}$ the usual basis of
$H^0(C,\Omega_1)$, and by $(\Omega_1,\Omega_2)$ a period matrix for
this basis with respect to some symplectic basis of $H_1(C,\Z)$, so
that $Z:=\Omega_1^{-1}\Omega_2\in\HH_g$. The Jacobian variety of $C$
can be described as the complex torus $J(C):=\C^g/(1_g|Z)$. We will
denote by $\Pi$ the normalized degree $g-1$  Abel-Jacobi map, $\Pi:~
C_{g-1}\longrightarrow~J(C)$, whose image  $\Pi(C_{g-1})$ is
precisely the divisor on $J(C)$ cut out by the Riemann theta
function $\q(Z;z)$.

The choice of the basis $\omega _{1},...,\omega _{g}$ of the space of holomorphic
differential forms on $C$ provides a \textit{canonical map} from $C$ to $\mathbb{P}^{g-1}=\mathbb{P}%
H^{0}(C,\Omega ^{1})^{*}$, given by:
$$
\begin{array}{rll}
\phi :C & \rightarrow & \mathbb{P}^{g-1} \\
P & \rightarrow & \phi (P)=(\omega _{1}(P),...,\omega _{g}(P)).
\end{array}
$$
\noindent Note that if the the differential forms
$\omega_1,\dots,\omega _{g}$ are  defined over the same number field
$K$ as the curve, then the canonical map is also defined over $K$.
The following result (which in fact is valid for a general curve)
relates the canonical images of certain divisors with their images
through the Abel-Jacobi map (cf. \cite{Guardia}):

\begin{proposition}
\label{hiperpla} Let $P_{1},...,P_{g-1}\in C(\bar{K})$ such that the
divisor $D=P_{1}+...+P_{g-1}$ satisfies $l(D)=1$. The equation:
$$
H_{D}(X_{1},...,X_{g}):=\left(
\begin{array}{lll}\displaystyle
\frac{\partial \theta }{\partial z_{1}}(\Pi (D)) ,& ...,
&\displaystyle \frac{\partial \theta }{\partial z_{g}}(\Pi (D))
\end{array}
\right) \Omega _{1}^{-1}\left(
\begin{array}{c}
X_{1} \\
\vdots  \\
X_{g}
\end{array}
\right) =0
$$
\noindent determines a hyperplane $H_{D}$ of $\mathbb{P}^{g-1}$,
which contains the divisor $\phi(D)$ on the curve $\phi(C)$.
\end{proposition}

\section{Symmetric normal forms for hyperelliptic curves}
\label{seccio-simetric}

The normalization of the roots of the polynomial $f(X)$ defining an hyperelliptic curve $Y^2=f(X)$ has been traditionally done following Rosenhain: one sends three of the roots of $f(X)$ to 0, 1 and $\infty$. This normalization has a number of advantages, but it could be not the most natural one. We introduce here a new normal model for hyperelliptic curves; the symmetries of this model allow the simplification of some common tasks related to hyperelliptic curves, as we will see later.

\subsection{Symmetric equations}

We assume that we are working over a field $K$ of characteristic different from 2, and denote by
$\mu_{4g}(\overline K)=\{\zeta_1,\dots,\zeta_{R}\}$ the ${4g}$-th roots of unity in $\overline K$.

\begin{definition} Let $C: Y^2=f(X)=(X-\alpha_1)\cdots(X-\alpha_{2g+2})$ be an hyperelliptic curve of genus $g$, defined over a field $K$ with $\operatorname{char} K\neq 2$.
For $i\neq j\in\{1,2,\dots,2g+2\}$ and $t\in\{1,\dots,R\}$ we define the following {\em symmetric invariants}:
\begin{itemize}
\item[a)] The {\em symmetric ratios } of $C$:
$$
p_{ijt}:=\zeta_t\sqrt[2g]{\prod_{k\neq i,j}\frac{\alpha_j-\alpha_k}{\alpha_i-\alpha_k}}=\zeta_t\sqrt[2g]{-\frac{f'(\alpha_j)}{f'(\alpha_i)}}\in\overline{K}.
$$
\item[b)] The {\em symmetric roots } of $C$:
$$
\ell_{ijtk}:=p_{ijt} \frac{\alpha_i-\alpha_k}{\alpha_j-\alpha_k}, \qquad k\in\{1,2,\dots,2g+2\}, k\neq i,j.
$$
\item[c)]
The {\em  symmetric normal models} for $C$:
$$
\begin{array}{rl}
\M_{ijt}: Y^2=\FF_{ijt}(X):=&X\prod_{k\neq i,j} (X-\ell_{ijtk})
\\\\
=&X^{2g+1}+G_{ijt,1}X^{2g}+\dots+G_{ijt,2g-1}X^2\pm X.
\end{array}
$$
(The coefficients $G_{ijt,k}$ will be called {\em symmetric coefficients}.)
\item[d)]
The {\em symmetric discriminants} of $C$:
\begin{equation}\label{modulardiscriminants}
\D_{ijt}=\prod_{r<s} (\ell_{ijtr}-\ell_{ijts})^2=\pm \frac{(\alpha_i-\alpha_j)^{2g(2g+1)}\Delta(f)}{f'(\alpha_i)^{2g+1}f'(\alpha_j)^{2g+1}}
\end{equation}
\end{itemize}
\end{definition}

\begin{remarks}$\,$
\begin{itemize}
\item[a)] The word {\em symmetric} refers to the relative position of the non-zero roots of $X(X^{2g}+G_1 X^{2g-1}+\dots+G_{2g-1}X\pm1)$ with respect to 0 and $\infty$.
\item[b)] In the case that the polynomial $f(X)$ defining the curve $C$ has degree $2g+1$, i.e., that one of its roots is $\alpha_i=\infty$, we can compute all the symmetric invariants with the same formulas, just substituting any factor $\alpha_i-\alpha_r$ by a 1.
\item[c)] Since for fixed $i,j,t$ the symmetric roots $\ell_{ijtk}$ are obtained from $\alpha_1, \dots, \alpha_{2g+2}$ by means of a common M\"{o}bius transformation, it is clear that $\M_{ijt}$ is a model of the curve $C$.

\item[d)] The roots of unity involved in the definition of the symmetric invariants are necessary to cover all the Galois conjugates of a given invariant. In order to simplify the notation, we will not write them explicitly anymore: we will denote the symmetric roots simply by $\ell_{ijk}$, assuming that a common root of unity has been chosen for fixed $i, j$. Hence, any equality involving the symmetric invariants should be understood, unless explicitly stated, modulo these $4g$-th roots of unity. For instance, when we write
$$
\ell_{ijk}=\ell_{jik}^{-1}
$$
it should be understood that for a proper choice of $\zeta_t,
\zeta_{t'}$ we have $\ell_{ijtk}=\ell_{jit'k}^{-1}$ for every  value
of $k$.
\end{itemize}
\end{remarks}

\begin{lemma} The symmetric roots satisfy the following relations:
\begin{itemize}
\item[a)]   $\ell_{jik}=\ell_{ijk}^{-1}$;
\item[b)] $\ell_{ijk}\ell_{jki}\ell_{kij}=-1$;
\item[c)] $\ds\ell_{irj}=\ell_{ijr}\prod_{k\neq i,j,r}(\ell_{ijk}-\ell_{ijr})$;
\item[d)] $\displaystyle \ell_{ijr}=\frac{\ell_{sji}-\ell_{sjr}}{\ell_{sij}}$;
\end{itemize}
\end{lemma}

As a consequence of part {\it a}), we see that
 when the symmetric model $\M_{ij}$ is
$$
\M_{ij}: Y^2=X^{2g+1}+G_1X^{2g}+\dots+G_{2g-1}X^2+ X,
$$
the symmetric model $\M_{ji}$ is
$$
\M_{ij}: Y^2=X^{2g+1}+G_{2g-1}X^{2g}+\dots+G_{1}X^2+X.
$$

There are also quite simple relations between  the   symmetric discriminants:

\begin{lemma} $\,$
\begin{itemize}
\item[a)]  $\D_{ij}=\D_{ji}$.\\
\item[b)] $\D_{ij}=\ell_{jki}^{2g(2g+1)}\D_{ik}$.\\
\item[c)] $\D_{ij}=\ell_{irs}^{2g(2g+1)}\ell_{jsi}^{2g(2g+1)}\D_{rs}$.
\end{itemize}
\end{lemma}

The symmetric normal model $\M_{ij}$ is determined by the choice of the roots $\alpha_i, \alpha_j$, which can be done in $(2g+2)(2g+1)/2$ different ways, and the choice of a $4g$-th root of unity, so that  we have up to $2g(2g+2)(2g+1)$  symmetric models for a generic hyperelliptic curve.
For arithmetic applications, it is worth noting that they may be not defined over the field of definition of the curve.

\begin{example} We have studied the symmetric models of elliptic curves in \cite{GuToVe}.
For an  elliptic curve  $E: Y^2=(X-e_1)(X-e_2)(X-e_3)$, the symmetric roots take the aspect:
$$
\ell_{ijr}=\sqrt{\pm\frac{e_i-e_r}{e_j-e_r}},
$$
and hence they are essentially the well-known {\em moduli} for $E$,  which are the roots of the equation
$$
256(k^4-k^2+1)^3-k^4(k^2-1)^2j_E=0,
$$
where $j_E$ is the absolute invariant of the elliptic curve $E$.
\end{example}

This fact generalizes to  hyperelliptic curves of any genus $g$: their symmetric roots  are absolute invariants of the curve  with certain level structure:

\begin{theorem}
\label{invariancia}
If
two hyperelliptic curves   defined over a field of odd characteristic are isomorphic, then their sets of symmetric roots are equal (after a proper labelling of the roots).
\end{theorem}

\noindent
{\bf Proof:} Suppose that we are given two isomorphic curves over a field $K$
$$
C: Y^2=\prod_i (X-\alpha_i), \qquad C':Y^2=\prod_i (X-\alpha'_i),
$$
with an isomorphism between them realized by a fractional linear transformation
 $\gamma(X)=\frac{AX+B}{CX+D}$  with  $A, B, C, D\in \overline{K}$ such that $\gamma(\alpha_j)=\alpha'_j$. We have:
$$
\frac{\alpha_i'-\alpha_k'}{\alpha_j'-\alpha_k'}=\frac{C \alpha_j+D}{C \alpha_i+D}\cdot
\frac{\alpha_i-\alpha_k}{\alpha_j-\alpha_k},
$$
and hence  the symmetric roots $\ell_{ijk}$ of $C$ and the symmetric roots  $\ell_{ijk}'$ of $C'$ will coincide. $\Box$

The symmetric roots $\ell_{ijk}$ being invariants of the curve $C$, any rational expression in them will produce new invariants. Particularly interesting will be the symmetric discriminants $\D_{ij}$.

\subsection{Reduction properties of symmetric models}$\quad$

We shall work now on a discrete valuation ring $A$, with field of
fractions $K$ of  characteristic different of 2. An integral model
for a hyperelliptic curve $C$ over $K$ can be given by an equation
of the form $Y^2=f(X)$ with $f(X)\in A[X]$. Such a model  can be
reduced modulo $\p$ (the prime ideal in $A$), yielding a new curve
$\tilde C$ over the residual field $k=A/\p$. This curve is
non-singular if and only if $\p\nmid 2\Delta(f)$, where $\Delta(f)$
denotes the discriminant of the polynomial $f(X)$; in this case it
is said that the curve $C$ has good reduction; otherwise it is said
that the curve has bad reduction. A minimal model for $C$ is an
integral model such that $\Delta(f)$ has minimal valuation with
respect to $\p$. A curve $C$ with bad reduction may have a model
over an extension $A'$ of $A$ with good reduction at the prime $\p'$
of $A'$ lying over $\p$; in this case it is said that $C$ has
potentially good reduction over $\p$.

The following results illustrate the interest of  symmetric models concerning the reduction of curves:

\begin{theorem}
Let $C: Y^2=f(X)$ be an hyperelliptic curve  over  $A$, with potentially good reduction. Let $A':=A[G_{ij,1},\dots, G_{ij,2g-1}]$ be the ring of definition of the symmetric model $\M_{ij}$ of $C$, and let $\p'$ be the prime in $A'$ lying over $\p$.
\begin{enumerate}
\item[a)] The symmetric coefficients $G_{ij}$ are $\p'$-integral.
\item[b)] The symmetric equation $\M_{ij}$  has good reduction at $\p'$.
\end{enumerate}
\end{theorem}
{\bf{Proof:}}
Let $\alpha_1,\dots,\alpha_{2g+2}$ be the roots of  $f(X)$ in $\overline K$.
Since
$$
\ell_{ijk}^{2g}=\pm\frac{(\alpha_i-\alpha_k)^{2g-1}}{(\alpha_j-\alpha_k)^{2g-1}}
{\prod_{r\neq i,j,k}\frac{\alpha_j-\alpha_r}{\alpha_i-\alpha_r}},
$$
it is clear that the symmetric roots $\ell_{ijt}$ are integral over the ring $A[\frac1{\Delta(f)}]$, since the denominators appearing in the last expression divide $\Delta(f)$. Thus the symmetric coefficients $G_{ij,k}$ are also integral over this ring.
Let $Y^2=f_1(X)=\prod_i(X-\beta_i)$ be a model of $C$ over a finite extension $A_1$ of $A$, with good reduction at the prime $\p_1\in \Spec(A_1)$ above $\p$. The discriminant $\Delta(f_1) $ must be a unit in $A_1$. We may now compute the symmetric models from this new model, since they are invariants of the curve $C$ by theorem \ref{invariancia}. We see thus that the coefficients $G_{ij,k}$ are integral over the ring $A_1[\frac1{\Delta(f_1)}]=A_1\supseteq A'$, and hence they are finally $\p'$-integral.

The  discriminant $\D_{ij}$ of the symmetric model $\M_{ij}$ is given by
$$
\label{Dij}
\begin{array}{rl}
\D_{ij,t}=& \ds p_{ij\,t}^{2g(2g-1)}(\beta_i-\beta_j)^{2g(2g-1)}\prod_{\begin{array}{c}
_{} ^{r<s}\\ {}^{r,s\neq i,j}
\end{array} }
\frac{(\beta_r-\beta_s)^2}{(\beta_j-\beta_r)^2(\beta_j-\beta_r)^2}
\\ \\
=& \ds \frac{(\beta_i-\beta_j)^{2g(2g-1)}
\ds\prod_{\begin{array}{c} ^{r\neq s}\\ {}^{r,s\neq i,j}\end{array}} (\beta_r-\beta_s)}
{\ds\prod_{r\neq i,j}(\beta_i-\beta_r)^{2g-1}(\beta_j-\beta_r)^{2g-1}}\in A_1[\frac1{\Delta(f_1)}]=A_1,
\end{array}
$$
so that it does not belong to $\p_1$, and hence the symmetric model $\M_{ij}$ has good reduction at $\p'$.
 $\Box$

\begin{corollary}
If  $v_\p(\D_{ij})\le 0$ for some $i, j$ then $C$ cannot have potentially good reduction at $\p$.
\end{corollary}

This corollary can be understood as a generalization of the well-known
  criterion for an elliptic curve having potentially good reduction   (\cite[p. 181]{Silverman}).
  It leads to the following definition:

\begin{definition}
The odd geometric locus of bad reduction of an hyperelliptic curve $C: Y^2=f(X)$ defined over a domain $A$ is
$$
\operatorname{BR}(C)^{odd}=\left\{\p\cap A  \mid\, \exists i,j  \mbox{ such that } \p\in\Spec\,A[\D_{ij}]
\mbox{ and }  v_\p(\D_{ij})<0 \right\}.
$$
\end{definition}

For the primes in $\operatorname{BR}(C)^{odd}$, symmetric models have also good properties:

\begin{theorem}
\label{minimal}
Let $C$ be a hyperelliptic curve over discrete valuation ring $A$, and suppose that
$C$ has not potentially good reduction at the unique  prime $\p$ in $\Spec A$.
If the symmetric coefficients $G_{ij,2},\dots, G_{ij,2g-1}$ are integral over $A$, then
the symmetric model $\M_{ij}$ is a minimal model for $C$
over the ring $A'=A[G_{ij,2},\dots, G_{ij,2g-1}]$.
\end{theorem}
{\bf{Proof:}}
Let $B=A'[\{\ell_{ij,r}\}_r]$, and let $\cP$ the prime of $B$ above $\p$.
By \cite[lemma 2.3]{Lockhart}, it is enough to see that the roots of the polynomial $X\prod_{r\ne i,j}(X-\ell_{ijr})$ defining $\M_{ij}$ are not all congruent  $\mod \cP$. This is clear, since $\prod_{r\ne i,j}\ell_{ijr}=\pm1$, so that these roots cannot be $0 \pmod \cP. \,\Box$

The minimality of symmetric models suggest that their coefficients $G_k$ should be small in some sense.
The following example illustrates this behavior:

{\bf Example:}
Weng (\cite{Weng}) computed the equation of a genus~3 hyperelliptic curve $C$ whose jacobian has complex multiplication by the field $K=\Q(w,i)$, where $w^3-w^2-10w+8=0, i^2=-1$. She found
$$
C:Y^2=f(X):=X^7 + 961X^5 - 3694084X^3 + 1832265664X.
$$
The discriminant of the polynomial $f(X)$ is $\Delta(f)=-2^{44}31^{35}$. A symmetric model for this equation is:
$$
Y^2=g(X):=X^7+\frac{\sqrt[3]{31}}4\,X^5  -\frac{\sqrt[3]{31^2}}4X^3  +X,
$$
which has only bad reduction at the primes dividing 2 in $\Q(\sqrt[3]{31})$, since $\Delta(g)=-2^{14}$.

\subsection{The $\mu$-invariants}
\label{mu-invariants}

We shall explain now a particular construction of the symmetric roots of a hyperelliptic curve, useful when we have not an explicit set of Weierstrass points, but certain intermediate invariants. This construction will be necessary in section \ref{Nullwerte-Weierstrass}.

\begin{definition} For $i,j,r,s,\in \{1,\dots,2g+2\}$
$$
\mu_{ijrs}:=\dfrac{(\alpha_i-\alpha_r)(\alpha_j-\alpha_s)}{(\alpha_i-\alpha_s)(\alpha_j-\alpha_r)}.\\
$$
\end{definition}

\begin{lemma}\quad
\label{lema-mu-invariants}
\begin{enumerate}[]
\item{a)} $\mu_{ijrs}=\mu_{rsij}=\mu_{jisr}=\mu_{jirs}^{-1}=\mu_{ijsr}^{-1}$.
\item{b)} $\ds\frac{\ell_{ijr}}{\ell_{ijs}}=\mu_{ijrs}$.
\item{c)} $\ell_{ijk}^{4g}=\prod_{r\neq i,j,k}\mu_{ijkr}^2$.
\end{enumerate}
\end{lemma}

\begin{proof}
Parts {\it a)} and {\it b)} are immediate. Part {\it b)} follows from:
$$
\begin{array}{rl}
\ell_{ijk}^{4g}&\displaystyle=\left(\dfrac{\alpha_i-\alpha_{k}}{\alpha_j-\alpha_{k}}\right)^{4g}\prod_{r\neq i,j,k}
\left(\dfrac{\alpha_j-\alpha_{r}}{\alpha_i-\alpha_{r}}\right)^2
\\\\
&=\displaystyle\prod_{r\neq i,j,k}
\left(\dfrac{(\alpha_i-\alpha_k)(\alpha_j-\alpha_s)}{(\alpha_i-\alpha_s)(\alpha_j-\alpha_k)}\right)^2=\prod_{r\neq i,j,k}\mu_{ijkr}^2.
\end{array}
$$
\end{proof}

Later on we will see how to compute the $\mu_{ijrs}$ by means of Jacobian Nullwerte. We see now how to obtain a symmetric equation from them.

\begin{proposition}\label{sym-ded}
A symmetric equation for an hyperelliptic curve can be deduced from a family $\{\mu_{ijrs}\}_{r,s}$.
\end{proposition}
\begin{proof}
We can compute the symmetric roots $\{\ell_{ijr}\}_r$ up to a $4g$-th root of unity by means of the third formula of the lemma. This root of unity must be determined coherently for all the symmetric roots, but this can be done using the second part of the lemma.
\end{proof}

\section{Preliminaries on Thetanullwerte and Jacobian Nullwerte}

\subsection{Theta-functions and theta-characteristics}

For $z\in\C^g, Z\in\HH_g$, the Riemann theta function on $\C^g$ is defined by:
$$
\theta (z,Z):=\sum_{n\in \mathbb{Z}^{g}}\exp (\pi i^{t}n.Z.n+2\pi i^{t}n.z).
$$
When the matrix $Z$ is fixed and we  consider $\q(z,Z)$ as a function of the first variable $z$, we will write it as $\q(z)$.
A vector $m={}^t\left(m', m''\right)$, with  $m',m''\in \R^g$, defines a translate of $\q$ as follows:
$$
\begin{array}{rl}
\q[m](z,Z):=& e^{\pi i ^tm'.Z.m'+2\pi i ^tm'.(z+m'')}\q(z+Zm'+m'') \\ \\
=&\ds\sum_{n\in \mathbb{Z}^{g}}e^{\pi i^{t}(n+m').Z.(n+m')+2\pi i^{t}(n+m').(z+m'')}.
\end{array}
$$
It is called theta function with characteristic $m$. If $m', m''\in \frac12\Z^g$ we call $m$ a {\em
 theta-charac\-te\-ris\-tic}. For a theta-charac\-te\-ris\-tic $m$, the corresponding theta function $\q[m](z,Z)$
 is an even or odd function of $z$ according to the parity of $m$, which is defined to be the parity of  the sign $e(m):=(-1)^{4^tm'.m''}$.

For a fixed matrix $Z\in \HH_g$, theta-characteristics are in bijection with two torsion points on the complex torus $T_Z:=\C^g/(1_g |Z)$:
$$
m=\left(\begin{array}{c} m' \\ m''\end{array}\right)\in\{0,1/2\}^{2g} \longleftrightarrow w:=Z.m'+m''\in T_Z[2].
$$
Along the paper, we shall use the symbols $m, m_k,\dots $ for theta-charac\-te\-ris\-tics, and the symbols $w, w_k,\dots$ for the corresponding 2-torsion points on $T_Z$, the relation between them being implicitly assumed. For instance, we define the parity of $w\in T_Z[2]$ to be the parity of $m$.

The values $\q[m](0,Z)$ are usually called {\em Thetanullwerte}, and
denoted shortly by $\q[m](Z)$ or $\q[w](Z)$ (even by  $\q[m]=\q[w] $
when  the matrix $Z$ is fixed). It is also usual to look at the
Thetanullwerte $\q[m](0,Z)$ as a functions of $Z$, i.e. defined on
the Siegel upper half space $\HH_g$. The $\C$-algebra spanned by
them is called {\em ring of Thetanullwerte}, and it is denoted by
$\C[\q]$. It has theoretical significance in relation to certain
rings of Siegel modular forms.

\subsection{Jacobian Nullwerte}

If we fix an odd characteristic $m$,  then the Thetanullwerte $\q[m](Z)$ vanishes  for every $Z$. For a sequence $M=\{m_1, \dots, m_g\}$ of  odd characteristics  one considers the matrix:
$$
J[M](Z):=J[m_1,\dots,m_g](Z):=\left(\begin{array}{ccc}
\displaystyle\frac{\partial \theta [m _{1}]}{\partial z_{1}}(0;Z) & \cdots  &
\displaystyle\frac{\partial \theta [m _{1}]}{\partial z_{g}}(0;Z) \\
\vdots  &  & \vdots  \\
\displaystyle\frac{\partial \theta [m _{g}]}{\partial z_{1}}(0;Z) & \cdots  &
\displaystyle\frac{\partial \theta [m _{g}]}{\partial z_{g}}(0;Z)
\end{array}
\right)
$$
and its determinant:
$$
[ m _{1},...,m _{g}](Z):=\pi^g D(M)(Z):= \det J[m_1,\dots,m_g](Z),
$$
which is usually called {\em Jacobian Nullwert}.
If the matrix $Z$ is fixed, we  shall denote $J[m_1,\dots,m_g](Z)$ and its determinant  by
$J[w_1,\dots,w_g]$ and $[ w _{1},...,w _{g}]$  respectively.

\subsection{Fundamental systems}

Given three theta characteristics $m_1$, $m_2$, $m_3$, define
$e(m _{1},m _{2},m _{3}) :=e(m _{1})e(m _{2})e(m
_{3})e(m _{1}+m _{2}+m _{3})$.
A tri\-plet $\{m _{1},m _{2},m _{3}\}$ is called {\em azygetic}  if $e(m _{1},m _{2},m _{3})=-1$, and syzygetic otherwise.  A sequence  $\{m_1,\dots,m_r\}$ is azygetic if every triplet contained in it is azygetic.

A \textit{fundamental system}  is an azygetic  sequence $S=\{m _{1},...,m _{2g+2}\}$ of $2g+2$
theta characteristics. A {\em special fundamental system} is a fundamental system with the first $g$ odd terms and the remaining $g+2$ even terms. The same concepts for two-torsion points on an abelian variety are defined analogously.

Fundamental systems play a basic role in the generalizations of
Jacobi's derivative formula obtained by Igusa. For low dimension we
have

\begin{theorem}[\cite{Igusa-Jacobi formula}]
\label{teorema d'Igusa} Assume $g\le 5$. Let $m_{1},..., m_{g}$ be
odd analytic theta characteristics such that the function
$[m_{1},..., m_{g}](Z)$ is not identically zero and is contained in
the ring of Thetanullwerte $\mathbb{C}[\theta ]$. Then $m_{1},...,
m_{g}$ can be completed to form a fundamental system, and:
$$
[m_{1},..., m_{g}](Z)=\pi ^{g}\sum_{\{m_{g+1},...,m_{2g+2}\}\in \mathcal{S}}\pm \prod_{i=g+1}^{2g+2}\theta [m _{i}](0;Z),
$$
\noindent
where $\mathcal{S}$ is the set of all $(g+2)$-tuples $\{m_{g+1},...,m_{2g+2}\}$ of even theta characteristics such that $\{m_{1},...,m_{g},$ $m_{g+1},...,m_{2g+2}\}$ form a fundamental system. If $Z$ is
the period matrix of a hyperelliptic curve, there is exactly one non-zero
term in the sum of the right hand side of the equality.
\end{theorem}

For higher dimensions $g>5$, Igusa provides a broader version of the
formula above, relating a sum of Jacobian Nullwerte with a certain
sum of products of Thetanullwerte
(\cite{Igusa-Poincare},\cite{Igusa-multip}).

Particular cases of Igusa's theorem are:
\begin{theorem}[Rosenhain's Formula, \cite{Rosenhain}, \cite{Guardia}]
\label{Rosenhain}
For any $Z\in\HH_2$ and any pair of odd characteristics $m_1, m_2$
\begin{equation}
[m_1,m_2]=\pm \pi^2 \prod_{
\begin{array}{c}
m\, {\mathrm {odd}}, \\ m\neq m_1,m_2
\end{array}
}\q[m_1+m_2-m],
\end{equation}
where the sign does not depend on $Z$.
\end{theorem}

\begin{theorem} [Frobenius, \cite{Frobenius}, \cite{Guardia}]
\label{Frobenius-g3} Let $C$ be a hyperelliptic curve of genus 3, with
Weierstrass points $W_{1},...,W_{8}$.
Let $w_{ij}:=\Pi(W_i+W_j)$, $w_{ijkr}:=\Pi(W_i+W_j+W_k-W_r)$.
The following equality holds for every
triplet $W_{i},W_{j},W_{k}$:
$$
[w_{ik}, w_{ij}, w_{jk}]=\pm\pi ^{3}\prod_{r\neq
i,j,k}\theta [w_{ijkr}],
$$
where the sign does not depend on $Z$.
\end{theorem}

\subsection{A fundamental system for Jacobians of hyperelliptic curves}

The following construction will be useful later:

\begin{proposition}[\cite{Guardia}]
\label{sistema fonamental hipereliptic}
Let  $W_{1},\dots,W_{g}$ be $g$ different Weierstrass points on a hyperelliptic curve $C$, and
denote $W_{g+1},\dots,W_{2g+2}$ the remaining Weierstrass points. Consider the divisors
$$
\begin{array}{rll}
D=&\sum_{i=1}^{g}W_{i},& \\
D_{i}=&D-W_{i}, & i=1,...,g, \\
D_{i}=&D+W_{i}-2W_{2g+2},\qquad &i=g+1,...,2g+2; \\
\end{array}
$$
The images $w_1=\Pi(D_1),\dots, w_{2g+2}=\Pi(D_{2g+2})$ of these divisors through the Abel-Jacobi map
form a
fundamental system of 2-tor\-sion points in $J(C)$.
\end{proposition}

\section{Jacobian Thomae's formula}
\label{Thomae-seccio}

After propositions \ref{hiperpla} and \ref{sistema fonamental hipereliptic}, if we fix $g$ Weierstrass points $W_1=(\alpha_1,0),\dots,W_g=(\alpha_g,0)$ on $C$
and define  $w_1=\Pi(\sum_{i\neq 1}W_i),\dots,w_g=\Pi(\sum_{i\neq g}W_i)$, the Jacobian Nullwerte $[w_1,\dots,w_g]$ has a clear geometric interpretation: the rows of this determinant are the equations
 of the $g$ hyperplanes spanned by $g-1$ of the fixed Weierstrass points. It is then natural to ask if there is some Jacobian version of Thomae's formula,
 connecting  $[w_1,\dots,w_g]$ with the discriminant of the polynomial $(X-\alpha_1)\cdots(X-\alpha_g)$. It was a beautiful surprise
 to find what seems to be a forgotten result in last page of Thomae's original paper,   which provides this connection:

\begin{theorem}[Thomae, {\cite[p. 222]{Thomae}}]
\label{Thomae-J}
Fix a partition $\{W_1,\dots, W_{g-1}\}\cup\{W_{g},\dots, W_{2g+2}\}$  of the set of Weierstrass points of $C$,  and consider the odd two-torsion point   $w_o=\Pi(W_1+\dots+W_{g-1})$. Define
$F_{w_o,1}(X)=\prod_{i=1}^{g-1}(X-\alpha_i)$,
$F_{w_o,2}(X)=f(X)/F_{w_o,1}(X)$. The following relation holds\footnote{Note that the  term  $\sqrt{\det{\Omega_1}}$ is misplaced in Thomae's paper.}:
$$
2(2\pi)^{g/2} \operatorname{grad}\theta[w_0](Z)=(\Delta(F_{w_o,1})\Delta(F_{w_o,2}))^{1/8}\sqrt{\det \Omega_1}S(\alpha_1,\dots,\alpha_{g-1}).\Omega_1,
$$
where
$S(\alpha_1,\dots,\alpha_{g-1})=((-1)^{g-1}\prod\alpha_i, \dots,
-\sum \alpha_i,1)$ is the row-vector formed by the coefficients of the polynomial $F_{w_o}(X)$.
\end{theorem}

We also recall the {\em standard} Thomae's formula (\cite{Thomae}, \cite{Mumford-Tata2}):
\begin{theorem}
\label{Thomae} Fix a partition $\{W_1,\dots,W_{g+1}\}\cup\{W'_{1},\dots,W'_{g+1}\}$ of the set of Weierstrass points on $C$, and consider the even two-torsion point $w_e=\Pi(W_1+\dots+W_{g+1}-2W)$ on $J(C)$. Define $G_{w_e,1}(X)=\prod_{i=1}^{g+1}(X-\alpha_i)$,
$G_{w_e,2}(X)=f(X)/G_{w_e,1}(X)$. The following relation holds:
$$
\q[w_e](Z)= (2\pi)^{-g/2}\sqrt{\det\Omega_1}\sqrt[8]{\Delta(G_{w_e,1})\Delta(G_{w_e,2})}.
$$
\end{theorem}

The combination of both Thomae formulas yields the following nice  results:

\begin{proposition}\label{Thomae-Jgeneral}
Let $w_o=\Pi(W_1+\dots+W_{g-1})$ and $w_e=\Pi(W_1'+\dots+W'_{g+1})$ . Then
$$
\frac2{\theta[w_e](Z)}\operatorname{grad}\theta[w_o](Z)=
\sqrt[8]{\frac{\Delta(F_{w_o,1})\Delta(F_{w_o,2})}
{\Delta(G_{w_e,1})\Delta(G_{w_e,2})}}
S(\alpha_1,\dots,\alpha_{g-1}).\Omega_1
$$
\end{proposition}

\begin{theorem}\label{Thomae-J2}
Write
$W_{g}=(\gamma_1,0), W_{g+1}=(\gamma_2,0)$, and denote by $ W_{k}=(\alpha_{k},0), k \neq0$ the remaining Weierstrass points on $C$. Define
$w_e=\Pi(W_1+\dots+W_{g+1}-2W_t)$, $G_{w_e,2}(X)=(X-\alpha_{g+2})\cdots(X-\alpha_{2g+2})$.
Then
$$
\frac2{\theta[w_e](Z)}\operatorname{grad}\theta[w_o](Z)=
\left(\frac{G_{w_e,2}(\gamma_1)G_{w_e,2}(\gamma_2)}{F_{w_o,1}(\gamma_1)F_{w_o,1}(\gamma_2)}\right)^{1/4}
 S(\alpha_1,\dots,\alpha_{g-1}).\Omega_1.
$$
\end{theorem}

\begin{proof}
We need only to simplify the expression inside the eighth root appearing in proposition \ref{Thomae-Jgeneral}.
We note that $G_{w_e,1}(X)=(X-\gamma_1)(X-\gamma_2)F_{w_0,1}(X)$ and $F_{w_o,2}(X)=(X-\gamma_1)(X-\gamma_2)G_{w_e,2}(X)$. Hence
$$
\frac{\Delta(F_{w_o,1})\Delta(F_{w_o,2})}{\Delta(G_{w_e,1})\Delta(G_{w_e,2})}=\displaystyle
\frac{\Delta(F_{w_o,1})\Delta(G_{w_e,2})(\gamma_1-\gamma_2)^2 G_{w_e,2}(\gamma_1)^2G_{w_e,2}(\gamma_2)^2}
{\Delta(F_{w_o,1})(\gamma_1-\gamma_2)^2 F_{{w_o,1}}(\gamma_1)^2F_{{w_o,1}}(\gamma_2)^2\Delta(G_{w_e,2})}.
$$
\end{proof}

\section{A remark on Igusa's theorem \ref{teorema d'Igusa}}
\label{seccio-remarks-igusa}

Let us consider theorem \ref{teorema d'Igusa} for hyperelliptic period matrices $Z$.
Fix a Jacobian Nullwerte   $[w_1,\dots,w_g](Z)$ with no identically zero row.
Igusa's theorem asserts that it can be represented as a product of $g+2$ Thetanullwerte $\q[w_1'](Z)\dots\q[w_{g+2}']$,
whenever $[w_1,\dots,w_g]$ is contained in the ring of Thetanullwerte $\C[\q]$. Igusa states a non-representability result in \cite[p. 93]{Igusa-odd}.
The results in previous section provide a simple way to generate new results in this direction.
After theorem \ref{Thomae-J}, we know that
$$
[w_1, \dots,w_g](Z)=2^{-g}(2\pi)^{-g^2/2}\Delta(\{w_r \}_r)
\left(\begin{array}{c}
S(\{\alpha_r^1\}_r)\\
\vdots\\
S(\{\alpha_r^g\}_r)
\end{array}
\right).\Omega_1
$$
where the term $\Delta(\{w_r\}_r)$ denotes a product of differences of the $x$-coordinates $\alpha_j^i$ of the Weierstrass points involved in the $w_i$.
If this expression equals a product of even Thetanullwerte, the standard Thomae's formula \ref{Thomae} implies that the determinant $\det(\left( S(\{\alpha_r^i\}_r) \right)_{i})$ must factor as a product of factors $(\alpha_j-\alpha_k)$. A formal computation may check this condition quite easily.

\section{Jacobi's formula revisited}
\label{seccio-jacobi-remarks}

The Jacobi triple product identity is usually written as
$$
\q_{1}(\tau)'=\pi\q_2(\tau)\q_3(\tau)\q_4(\tau),
$$
where  $q_r(\tau)$ is the usual notation for the Thetanullwerte in dimension 1.
A proper rearrangement of the formula drives to interesting remarks.
Consider the elliptic curve  $E_\tau$  associated to the complex torus $\C/\langle 1,\tau\rangle$:
$$ E_\tau: Y^2=X^3-g_2(\tau)X-g_3(\tau)=(X-e_1(\tau))(X-e_2(\tau))(X-e_3(\tau)),
$$
 where $g_2(\tau), g_3(\tau)$ are given by the classical Eisenstein series, and $e_j(\tau)=\dfrac{\pi^2}3(\q_r[\tau]\pm\q_s[\tau])$ (cf. \cite[p. 133 ]{McKean-Moll}).
We have seen in  \cite{Guardia-MRL} that the differential form ($r,s\neq 1$)
$$
\frac{\q_1(\tau)'}{\q_t(\tau)}dz=\pi\q_r(\tau)\q_s(\tau)dz
$$
is  defined over a finite extension of $\Q(j(\tau))$, and it is well-known (cf. \cite[p. 132]{McKean-Moll}) that
$$
\pi\q_r(\tau)\q_s(\tau)=\sqrt[4]{(e_i(\tau)-e_j(\tau))(e_i(\tau)-e_k(\tau))}.
$$
Hence, we can re-write Jacobi formula as
$$
\frac{\q_1(\tau)'}{\pi\q_r(\tau)}=\sqrt[4]{(e_i(\tau)-e_j(\tau))(e_i(\tau)-e_k(\tau))}\in \overline{\Q(j(\tau))}.
$$
For a general elliptic curve $E: Y^2=(X-\alpha_1)(X-\alpha_2)(X-\alpha_3)$, we obtain
$$
\frac{\q_1(\tau)'}{\pi\q_r(\tau)}=\omega_1\sqrt[4]{(\alpha_i-\alpha_j)(\alpha_i-\alpha_k)},
$$
where $\omega_1$ is a proper period of $E$.

In this way, Jacobi formula can be thought as a result in the area of
algebraic values of transcendental functions. In this direction, proposition \ref{Thomae-Jgeneral} provides a general version of it:

\begin{theorem}
Let $C: Y^2=(X-\alpha_1)\cdots(X-\alpha_{2g+2})$ be a genus $g$ hyperelliptic curve, and let
$(\Omega_1|\Omega_2)$ a period matrix of $C$ with respect to a symplectic basis of $H_1(C,\Z)$, so that
$Z:=\Omega_1^{-1}\Omega_2\in\HH_g$. For any choice of $g$ odd two torsion points $w_1,\dots,w_g\in J(C)$ and $g$ even
two torsion points $w_1',\dots,w_g'\in J(C)$:
$$
\left(\frac{1}{\pi^g\det\Omega_1}\dfrac{[w_1,\dots,w_g]}{\q[w_{1}']\cdots\, \q[w_{g}']}\right)^8\in \Q(\alpha_1,\dots,\alpha_{2g+2}),
$$
\end{theorem}

\section{Algebraic differential forms and periods}
\label{matrius-algebraiques}

\subsection{Algebraic differential forms for abelian varieties}$\quad$
\label{diferencials-va}

Given a principally polarized abelian variety $A$ of dimension $g$, defined over a field $k\subset \C$, and a basis $\w_1,\dots,\w_g\in H^0(A,\Omega^1_{/k})$, we form the period matrix $(\Omega_1| \Omega_2)$ of these differential forms with respect to any symplectic basis of $H_1(A,\Z)$. We may identify $A(\C)$ with the complex torus $\C^g/(\Omega_1 | \Omega_2)$. It is well-known that the $4^g$ theta-functions $\q[m](z)$ with half integer characteristic provide a
map
$$
\begin{array}{rcl}
A(\C) & \longrightarrow &\P^{4^g-1}(\C)\\
z & \longrightarrow & (\q[m](\Omega_1^{-1}.2z))_m,
\end{array}
$$
whose image is isomorphic to $A$ over a certain field $K$ which is a finite extension of $k$. We now take an even theta function $\q_0(z)$ such that $\q_0(0)\neq0$, and an odd one $\q_1(z)$, and form the quotient $\q_1(z)/\q_0(z)$; it can be seen as an element of $K(A)$, the field of $K$-rational functions on $A$. Hence, its differential at the point $z=0$ must be a $K$-linear combination of the original differential forms $\w_1,\dots,\w_g$:
$$
d(\q_1(z)/\q_0(z))\mid_{z=0}=\q_0(0)^{-1}d\q_1(z)\mid_{z=0}=\sum_{k=1}^g c_k \w_k.
$$
We observe that it is always possible to choose $g$ odd theta functions $\q_1(z),\dots,\q_g(z)$ such that
$d\q_1(z),\dots, d\q_g(z)$ are linearly independent differential forms on the torus  $\C^g/(\Omega_1 | \Omega_2)$
(cf. \cite{Shimura-CM}, p. 192). These differential forms will be defined over $K$, but not necessarily over $k$. We will use the term {\bf algebraic differential form} to describe a differential form on a variety defined over a finite extension of its field of moduli.

It is a frequent setting that the only available data of a complex abelian variety $A$ is a normalized period matrix
$Z\in\HH_g$. In this situation,
in order to build a period lattice for $A$ coming from and algebraic basis of differential forms, we look for a set of theta-functions $\q_0(z), \q_1(z),\dots,\q_g(z)$ as before. Equivalently, we look for 2-torsion points $w_0, w_1,\dots, w_g$ such that
$$
\Omega_1(w_1,\dots, w_g; w_0; Z):=\frac1{2\pi i \q[w_0]} J[w_1,\dots,w_g](Z)
$$
is a non-singular matrix. By the remarks above, then
\begin{equation}\label{differentialforms}
\left(\widetilde\w_1,\dots,\widetilde\w_g \right)=(dz_1,\dots,dz_g) \Omega_1(w_1,\dots, w_g; w_0;Z)
\end{equation}
are algebraic differential forms and
yield a basis of $H^0(A,\Omega^1_{/\overline{K}})$, with periods  $\Omega_1(w_1,\dots, w_g; w_0;Z)(1_g \mid Z)$.

\subsection{Algebraic differential forms for hyperelliptic curves}$\quad$

We assume now that $A$ is the jacobian variety $J(C)$ of a hyperelliptic curve $C: Y^2=f(X)$ of genus $g$, defined over a field $k\subseteq \C$. We have studied how to find a good algebraic equation for an elliptic curve from its normalized period lattice in \cite{Guardia-MRL}, so that from now on we will assume that we are working with a hyperelliptic curve of genus $g\ge2$.

Our initial data will be a normalized period matrix $Z\in\HH_g$ for the curve $C$, coming from a  model $Y^2=(X-\alpha_1)\cdots(X-\alpha_{2g+2})$ of the curve $C$.
    The  procedure described above provides a basis $\left(\widetilde\w_1,\dots,\widetilde\w_g \right)=\Pi^\ast\left((dz_1,\dots,dz_g) \Omega_1(w_1,\dots, w_g; w_0;Z)\right)$ of algebraic differential forms  in $H^0(C,\Omega^1_{/\overline{K}})$ derived from a set of $g$ odd two-torsion points $w_1,\dots,w_g \in J(C)$, which  are only subject to the condition $[w_1,\dots,w_g](Z)\neq0$.

We can easily describe now a geometric method to build a basis $\widetilde\w_1,\dots,\widetilde\w_g$ of algebraic holomorphic differentials from a normalized period matrix for $C$.
We take $g$ Weierstrass points  $W_1=(\alpha_1,0),\dots,W_g=(\alpha_g,0)$  on $C$, and form the divisors $D_i=\sum_{j\ne i }W_j$, and their images $w_i=\Pi(D_i)$ in $J(C)$ through the Abel-Jacobi map.
By \ref{hiperpla}, the rows of the matrix $J[w_1,\dots,w_g]$ are linearly independent, and thus for every even theta-characteristic $\q[w_0]\not\equiv0$
$$
\left(\widetilde\w_1,\dots,\widetilde\w_g \right)=\Pi^\ast((dz_1,\dots,dz_g) \Omega_1(w_1,\dots, w_g; w_0;Z))
$$
is a basis of $H^0(C,\Omega^1_{/\overline{K}})$. The point is  that this construction can be done
working only with the two torsion points of $J(C)$,
without explicit knowledge of the Weierstrass points of $C$. In this situation, theorem \ref{Thomae-Jgeneral} will determine  the field of definition of these differential forms.

\section{Jacobian Nullwerte and symmetric equations}$\quad$
\label{Nullwerte-Weierstrass}

 We now describe  a method to build a symmetric equation of $C$ by means of Jacobian Nullwerte. The fundamentals of the method are certain relations between quotients of Jacobian Nullwerte and the $\mu$-invariants. Although these relations can be deduced directly from theorem \ref{Thomae-J} or theorem \ref{Thomae-Jgeneral}, we provide some geometric intuition to derive these formulas.

We assume that the period matrix $Z$ comes from a certain model $Y^2=(X-\alpha_1)\cdots(X-\alpha_{2g+2})$ of the curve $C$.

    The  procedure described in section \ref{matrius-algebraiques} provides a basis $\left(\widetilde\w_1,\dots,\widetilde\w_g \right)=\Pi^\ast\left((dz_1,\dots,dz_g) \Omega_1(w_1,\dots, w_g; w_0;Z)\right)$ of algebraic differential forms  in $H^0(C,\Omega^1_{\overline{K}})$ coming from certain odd two-torsion points $w_1,\dots,w_g \in J(C)$, which  are only subject to the condition $[w_1,\dots,w_g](Z)\neq0$.
 By the Riemann singularity theorem, we know that  $w_i=\Pi(D_i)$ for certain geometric theta-characteristic $D_i$ on $C$ with  $l(D_i)=1$. This theta-characteristic  $D_i$ must be the sum of $g-1$ Weierstrass points (cf. \cite[p. 288]{ACGH84}).

\begin{proposition}
\label{imatge-canonica}
Let  $W$ be a Weierstrass point on $C$.
The image of $W$ through the canonical map $\phi_G$ given by $\widetilde\w_1,\dots,\widetilde\w_g$ is
$$\phi_G(W)=\left(
[w_1, w_2',\dots, w_g']:
[w_2, w_2',\dots, w_g']:\dots:
[w_g, w_2',\dots, w_g']
\right),
$$
where $w'_j=\Pi(W+\sum_{r=2}^{g-1} W_{jr})$ are   any $g-1$ odd 2-torsion points on $J(C)$ whose associated divisors contain the point $W$ and which do not coincide with the $w_i$.
\end{proposition}

\begin{proof}
Let us write $\Omega_1=\Omega_1(w_1,\dots,w_g,w_0;Z)$. By \cite[proposition 3.1]{Guardia}, the solution of the linear system
$$
\left(
\begin{array}{ccc}
\frac{\partial \theta }{\partial z_{1}}(w_{2}',Z) & \cdots  &
\frac{\partial
\theta }{\partial z_{g}}(w_{2}',Z) \\
\vdots  &  & \vdots  \\
\frac{\partial \theta }{\partial z_{1}}(w_{g}',Z) & \cdots  &
\frac{\partial
\theta }{\partial z_{g}}(w_{g}',Z)
\end{array}
\right) \Omega _{1}^{-1}\left(
\begin{array}{c}
X_{1} \\
X_{2} \\
\vdots  \\
X_{g}
\end{array}
\right) =\left(
\begin{array}{c}
0 \\
0 \\
\vdots  \\
0
\end{array}
\right),
$$
  is the canonical image $\phi_G(W)$ of the point $W=W_1'$. Let $(Y_1,\dots,Y_g)=(X_1,\dots,X_g) {}^t (\Omega_1 {^t})^{-1}$, and call $A$ the first matrix in the equality above. The solutions of the system
$$
A\left(
\begin{array}{c}
Y_{1} \\
\vdots  \\
Y_{g}
\end{array}
\right) =\left(
\begin{array}{c}
0 \\
\vdots  \\
0
\end{array}
\right)
$$
are:
$$
(Y_1:Y_2:\dots:Y_g)=(A_1:-A_2:\dots:(-1)^{g+1}A_g),
$$
where $A_j$ is the determinant of the matrix obtained by deleting the $i$-th column of $A$.
Now:
$$
\left(
\begin{array}{c}
X_{1} \\
\vdots  \\
X_{g}
\end{array}\right)=
\left(
\begin{array}{ccc}
\frac{\partial \theta }{\partial z_{1}}(w_{1},Z) & \cdots  &
\frac{\partial
\theta }{\partial z_{g}}(w_{1},Z) \\
\vdots  &  & \vdots  \\
\frac{\partial \theta }{\partial z_{1}}(w_{g},Z) & \cdots  &
\frac{\partial
\theta }{\partial z_{g}}(w_{g},Z)
\end{array}
\right)
\left(
\begin{array}{c}
A_{1} \\
\vdots  \\
(-1)^{g+1}A_{g}
\end{array}\right),
$$
and the result follows immediately.
\end{proof}

This result is specially significant when  the original  theta-charac\-te\-ris\-tics $D_1,\dots, D_g$ are well-posed:

\begin{proposition}
\label{mus}
Suppose that the 2-torsion points $w_1,\dots,w_g$ related to the basis $\widetilde\w_1,\dots,\widetilde\w_g$ of $H^0(C,\Omega^1_{\bar K})$ are the images of certain divisors $D_i=\sum_{j\neq i}W_j$.
Let $\phi_G(W_r)=(u_1,\dots,u_g), \phi_G(W_s)=(v_1,\dots,v_g)$, with $r,s>g$. Then
$$
\frac{u_m v_n}{u_n v_m}=\mu_{mnrs}.
$$
\end{proposition}

\begin{proof}
Let us denote by $\phi_S: C\rightarrow \P^{g-1}$ the canonical map given by the standard basis $\{\frac{x^rdx}{y}\}_r$ of $H^0(C,\Omega^1)$. We know that $\phi_S(W_n)=(1,\alpha_n,\dots\alpha_n^{g-1})$, and proposition \ref{hiperpla} shows that $\phi_G=P^{-1}\circ\phi_S$, where $P$ is the linear map given by the Vandermonde matrix:
$$
V(\alpha_1,\dots \alpha_g):=\left(%
\begin{array}{cccc}
    1 & 1 & \dots & 1 \\
  \alpha_1 & \alpha_2 & \dots & \alpha_{g} \\
  \alpha_1^2 & \alpha_2^2 &  & \alpha_{g}^2 \\
  \vdots & \vdots &  & \vdots \\
  \alpha_1^{g-1} & \alpha_2^{g-1} &\dots & \alpha_{g}^{g-1}
\end{array}%
\right).
$$
Let $F(X)=(X-\alpha_1)\cdots(X-\alpha_g)$, and consider the
polynomials
$H_i(X)=\frac1{F'(\alpha_i)}\frac{F(X)}{(X-\alpha_i)}=\sum_{j=0}^{g-1}a_{i,j+1}X^j$,
whose coefficients are  the entries of the matrix $IV(\alpha_1,\dots
\alpha_g):=V(\alpha_1,\dots,\alpha_g)^{-1}=\left(a_{i,j}\right)_{i,j=1\dots,g}$,
since $H_i(\alpha_j)=\delta_{ij}$. The map $P^{-1}$ is then given by
a matrix \linebreak
 $\mbox{diag}(\lambda_1,\dots,\lambda_g)IV(\alpha_1,\dots
\alpha_g)$, where the $\lambda_i$ are certain numbers which we don't
need to determine. Hence the coordinates of $\phi_G(W_r)$ are:
$$\phi_G(W_r)=
\left(\begin{array}{c}
\frac{\lambda_1 F(\alpha_r)}{(\alpha_r-\alpha_1)F'(\alpha_1)}\\ \vdots \\ \vdots\\ \frac{\lambda_g F(\alpha_r)}{(\alpha_r-\alpha_g)F'(\alpha_g)}
\end{array}\right),
$$
and the result follows immediately.
\end{proof}

We can combine this two results in order to find nice formulas for double ratios of Jacobian Nullwerte. An easy example is the following:

\begin{corollary}
\label{J-alfa}
Let $W_1,\dots,W_{g}, W_r, W_s$ be  $g+2$ different Weierstrass points on $C$, and form the divisors:
$$
\begin{array}{ll}
D=\sum_{i=1}^g W_i, &\\
D'=\sum_{i>g, i\neq r,s} W_i,&\\\
D_i=D-W_i, & i=1,\dots,g,\\
D'_i=D_i+W_{r}-W_1, & i=2,\dots,g,\\
D''_i=D_i+W_{s}-W_1, & i=2,\dots,g.\\
\end{array}
$$
\begin{itemize}
\item[a)] Let $w_i=\Pi(D_i)$, $w'_i=\Pi(D'_i)$, $w''_i=\Pi(D''_i)$. The following equalities hold:
\begin{equation}\label{quocientJ}
\frac{[w_m,w_2',\dots,w_g'][w_n,w_2'',\dots,w_g'']}{[w_m,w_2'',\dots,w_g''][w_n,w_2',\dots,w_g']}=\mu_{mnrs}.
\end{equation}
\item[b)] Let $w_{mr}=\Pi(D'+W_m+W_r)$, $w_{ms}=\Pi(D'+W_m+W_s)$,
$w_{nr}=\Pi(D'+W_n+W_r)$, $w_{mr}=\Pi(D'+W_n+W_s)$. We have:
$$
\frac{[w_m,w_2',\dots,w_g'][w_n,w_2'',\dots,w_g'']}{[w_m,w_2'',\dots,w_g''][w_n,w_2',\dots,w_g']}=\pm\left(
\frac{\q[w_{mr}]\q[w_{ns}]}{\q[w_{ms}]\q[w_{nr}]}\right)^2.
$$
\end{itemize}
\end{corollary}

\begin{proof}
The first equality follows from propositions \ref{imatge-canonica} and \ref{mus}. The second relation is derived from the first one and Thomae's formula \ref{Thomae}.
\end{proof}

As a by product of this corollary we obtain a relation between some double ratios of Jacobian Nullwerte:
\begin{corollary}
Consider a second family  $W'''_1,\dots, W'''_g$ of $g$ Weierstrass points  different from $W_{g+1},\dots,W_{g+2}$, and form the divisors
$D'''_i=\sum_{j\neq i} W'''_j$. Let $w'''_i=\Pi(D'''_i)$. We have:
$$
\frac{[w_m,w_2',\dots,w_g'][w_n,w_2'',\dots,w_g'']}{[w_m,w_2'',\dots,w_g''][w_n,w_2',\dots,w_g']}=
\frac{[w_m''',w_2',\dots,w_g'][w_n''',w_2'',\dots,w_g'']}{[w_m''',w_2'',\dots,w_g''][w_n''',w_2',\dots,w_g']}
$$
\end{corollary}

All the formulas above concerning Jacobian Nullwerte have been obtained by geometric means, but they could have been obtained directly from theorem \ref{Thomae-J}, which is a source for lots of such relations.

An important consequence of corollary \ref{J-alfa} and proposition \ref{sym-ded} is the fact that{\em  Jacobian Nullwerte provide an effective solution to the hyperelliptic Schottky problem in every genus}.
In the next two sections we shall explicit this construction for genus 2 and genus 3  hyperelliptic curves. These are the cases which are usually considered for applications (for instance, in cryptography) and which are at the reach of present standard computational power.

\section{Genus 2 curves}

We shall now consider
the results in previous sections for the particular case of genus 2 curves, where a number of refinements can be obtained, both in the algebraic and the analytic sides. We will explain how to find a symmetric equation for a hyperelliptic curve given  its Igusa-Clebsch  invariants or a normalized period matrix. Both methods are simple and efficient.

\subsection{Symmetric invariants for genus 2 curves}

Given a symmetric equation $Y^2=X(X^4+G_1X^3+G_2X^2+G_3X+1)$, it is a simple matter the determination of its Igusa-Clebsch invariants:
\begin{equation}
\label{GtoIgusa}
\begin{array}{rl}
I_2=& 2\,\left( 20 + 3\,G_2^2 - 8\,G_1\,G_3 \right),\\\\
I_4=&  -4\,\left( 20 + 3\,G_1^2\,G_2 - 9\,G_2^2 - G_1\,G_3 -
     G_1^2\,G_3^2 + 3\,G_2\,G_3^2 \right), \\\\
I_6=&  -2\,\left( 160 + 18\,G_1^4 - 13\,G_1^2\,G_2 - 88\,G_2^2 +
     12\,G_1^2\,G_2^3 - 36\,G_2^4 \right.\\
   &   - 32\,G_1\,G_3 -
     38\,G_1^3\,G_2\,G_3 + 119\,G_1\,G_2^2\,G_3 -
     14\,G_1^2\,G_3^2\\
     & - 13\,G_2\,G_3^2 -
     4\,G_1^2\,G_2^2\,G_3^2 + 12\,G_2^3\,G_3^2 +
     12\,G_1^3\,G_3^3 \\
     &\left. - 38\,G_1\,G_2\,G_3^3 +
     18\,G_3^4 \right),\\\\
I_{10}=&   - 27\,G_1^4 + 144\,G_1^2\,G_2 - 128\,G_2^2 -
     4\,G_1^2\,G_2^3 + 16\,G_2^4\\
     & - 192\,G_1\,G_3 +
     18\,G_1^3\,G_2\,G_3 - 80\,G_1\,G_2^2\,G_3 -
     6\,G_1^2\,G_3^2 \\
     &+ 144\,G_2\,G_3^2 +
     G_1^2\,G_2^2\,G_3^2 - 4\,G_2^3\,G_3^2 -
     4\,G_1^3\,G_3^3 \\
     &+ 18\,G_1\,G_2\,G_3^3 -
     27\,G_3^4+256.
\end{array}
\end{equation}

On the other hand, the determination of an hyperelliptic curve with prescribed invariants $I_2=\I_2, I_4=\I_4, I_6=\I_6, I_{10}=\I_{10}$ is a non-trivial problem, solved by \cite{Mestre} and \cite{Biel}. We will explain here an elementary method to find a symmetric equation with prescribed invariants, which takes profit of the simplicity of the expressions above.
Since the Igusa-Clebsch invariants are homogeneous invariants and the symmetric coefficients are absolute invariants, we need to introduce a proportionality constant, and solve the equations above for $I_{2k}=r^k \I_{2k}$.

First of all, we note that the formulas above are symmetric polynomials in  $G_1, G_3$, so that we can express them in terms of
\begin{equation}\label{G1G3S1S2}
S_1=G_1+G_3, \qquad S_2=G_1 G_3.
\end{equation}
We obtain:
$$
\begin{array}{rl}
r\I_2=&2\,\left( 20 + 3\,G_2^2 - 8\,S_2 \right),\\\\
r^2\I_4=&4\,\left( -20 + 9\,G_2^2 - 3\,G_2\,S_1^2 + S_2 +
    6\,G_2\,S_2 + S_2^2 \right),\\\\
r^3\I_6=&2\,\left( -160 + 88\,G_2^2 + 36\,G_2^4 + 13\,G_2\,S_1^2 -
    12\,G_2^3\,S_1^2 - 18\,S_1^4 \right.\\
    & + 32\,S_2 -
    26\,G_2\,S_2 - 119\,G_2^2\,S_2 + 24\,G_2^3\,S_2
    +     72\,S_1^2\,S_2 \\
    &+38\,G_2\,S_1^2\,S_2 \left.-
    22\,S_2^2 - 76\,G_2\,S_2^2 + 4\,G_2^2\,S_2^2 -
    12\,S_2^3 \right),\\\\
r^5\I_{10}=& - 128\,G_2^2 + 16\,G_2^4 + 144\,G_2\,S_1^2 -
  4\,G_2^3\,S_1^2 - 27\,S_1^4 - 192\,S_2 \\
  &-
  288\,G_2\,S_2 - 80\,G_2^2\,S_2 + 8\,G_2^3\,S_2 +
  108\,S_1^2\,S_2 + 18\,G_2\,S_1^2\,S_2 \\
  &-
  60\,S_2^2 - 36\,G_2\,S_2^2 + G_2^2\,S_2^2 -
  4\,S_2^3    +256.
\end{array}
$$
From the first equality we have:
\begin{equation} \label{S2}
S_2=(6G_2^2+40-r\I_2)/16,
\end{equation}
and we replace this relation in the remaining equations:
$$
\begin{array}{rl}
2^6\I_4=&36\,G_2^4 +576\,G_2^3 - 12\left(r \I_2-240 \right) \,G_2^2 -
  96\,\left(  8\, S_1^2+ r\I_2 -40 \right)\,G_2\\
  &+
  \left(r \I_2 -120\right) \,\left(r \I_2+24 \right)  ,\\\\

2^9\I_6=& - 72\,G_2^6 - 1728\,G_2^5  +
  12\,\left(  11\,\I_2\,r-1440  \right) \,G_2^4
  \\
  &+  192\,\left( 12\,S_1^2+ 11\,\I_2\,r   -492\right)\,G_2^3\\
  &-
  2\,\left(  - 13824\,S_1^2  + 19\,{\I_2}^2\,r^2 - 5856\,\I_2\,r+152640 \right)\,G_2^2 \\
  &-
  16\,\left( 19\,\I_2\,r-864 \right) \,\left( 8\,S_1^2 + \I_2\,r -40 \right)  \,{G_2}\\
&- 18432\,S_1^4 -
  4608\,\left( + \I_2\,r -40\right) \,S_1^2 + 3\,\I_2^3\,r^3
  \\
&- 448\,\I_2^2\,r^2 + 19392\,\I_2\,r -414720 ,
\\\\
2^{10}\I_{10}=&
 - 72\,G_2^6- 2112\,G_2^5+
  4\,\left(  15\,\I_2\,r-6344 \right) \,G_2^4  \\
 & +
  64\,\left( 44\,S_1^2 + 19\,\I_2\,r -2488 \right)\,G_2^3\\
 & - 2\,\left( - 20736\,S_1^2 + 7\,{\I_2}^2\,r^2- 4560\,\I_2\,r + 273600 \right)\,G_2^2  \\
&-   144\,\left( \I_2\,r -168\right) \,\left( 8\,S_1^2 + \I_2\,r -40 \right)  \,{G_2}\\
 & - 27648\,S_1^4 -   6912\,\left( + \I_2\,r -40\right) \,S_1^2\\
&+ {\I_2}^3\,r^3- 360\,{\I_2}^2\,r^2  + 36288\,\I_2\,r    -677376.
  \end{array}
$$
The first equality yields:
\begin{equation}\label{S12}
\begin{array}{l}
S_1^2=\frac{   36\,G_2^4 + 576\,G_2^3  - 12\,\left( r\,{\I_2} -240  \right)\,G_2^2  -
  96\,\left(  r\,{\I_2}-40 \right)\,{G_2}  +
r^2\,\I_2^2 -96\,{\I_2}- 64\,r^2\,{\I_4}
-2880 }{768G_2}.
\end{array}
\end{equation}

We now substitute this value in the equations for $\I_6, \I_{10}$:
\begin{equation}
\label{G2r}
\begin{array}{l}
_{
432\,G_2^8
-   864\,G_2^6\,\left( 48 + r\,{\I_2} \right)
+ 72\,G_2^4\,
   \left( 18240 + 672\,r\,{\I_2} + 5\,r^2\,\I_2^2 + 64\,r^2\,{\I_4}
     \right)
   }
   \\ _{
     - 8\,G_2^2\,\left( 1797120 + 81216\,r\,{\I_2} + 432\,r^2\,\I_2^2 +
     7\,r^3\,\I_2^3 + 27648\,r^2\,{\I_4} +
     1856\,r^3\,{\I_2}\,{\I_4} - 6144\,r^3\,{\I_6} \right)
     }
     \\
     _{
+2^{12}3^45\,r\,{\I_2} - 576\,r^3\,{\I_2}\,
   \left( \I_2^2 - 64\,{\I_4} \right)  +
  3\,r^4\,{\left( \I_2^2 - 64\,{\I_4} \right) }^2 +
  3456\,r^2\,\left( 3\,\I_2^2 + 320\,{\I_4} \right)+2^{12}3^55^2
=0,}
     \\ \\
     _{
  144\,G_2^8 -
  96\,G_2^6\,\left( 176 + 5\,r\,{\I_2} \right)
  +
  8\,G_2^4\,\left( 89280 + 4704\,r\,{\I_2} + 59\,r^2\,\I_2^2 +
     448\,r^2\,{\I_4} \right)
     }\\_{
       + 24\,G_2^2\,
   \left( -525312 - 36288\,r\,{\I_2} - 720\,r^2\,\I_2^2 -
     5\,r^3\,\I_2^3 - 9216\,r^2\,{\I_4} -
     192\,r^3\,{\I_2}\,{\I_4} + 8192\,r^5\,{\I_{10}} \right)
     }
     \\_{
+ 2^{12}3^5 5\,r\,{\I_2} - 1728\,r^3\,{\I_2}\,
   \left( \I_2^2 - 64\,{\I_4} \right)  +
  9\,r^4\,{\left( \I_2^2 - 64\,{\I_4} \right) }^2 +
  10368\,r^2\,\left( 3\,\I_2^2 + 320\,{\I_4} \right) +2^{12} 3^6 5^2=0.
   }
\end{array}
\end{equation}

We can solve these equations by a resultant elimination. It turns out that the elimination of $G_2$ produces a much simpler result:
\begin{equation}\label{r}
\begin{array}{l} _{2^83^6\,r^{15}\,\I_{10}^4 +
  2^63^6\,\,r^{13} \, \I_{10}^3\left( {\I_2}\,{\I_4} - 4\,{\I_6} \right)
  -  2^63^5 \,r^{12}\,\I_{10}^3 \,\left( \I_2^2 - 16\,{\I_4} \right)
   }\\
_{   +
  108\,r^{11}\,\I_{10}^2\,\left( 19\,\I_2^2\,\I_4^2 +
     8\,\I_4^3 - 168\,{\I_2}\,{\I_4}\,{\I_6} +
     360\,\I_6^2 + 5616\,{\I_2}\,{\I_{10}} \right)
     }
     \\
_{
 - 216\,r^{10}\,\I_{10}^2  \,
   \left( 11\,\I_2^3\,{\I_4} + 16\,{\I_2}\,\I_4^2 -
     36\,\I_2^2\,{\I_6} - 192\,{\I_4}\,{\I_6} -
     105408\,{\I_{10}} \right)
     }\\
_{+
  2\,r^9\,{\I_{10}}\,( \I_2^5\,\I_4^2 +
     25\,\I_2^3\,\I_4^3 - 26\,{\I_2}\,\I_4^4 -
     6\,\I_2^4\,{\I_4}\,{\I_6} -
     324\,\I_2^2\,\I_4^2\,{\I_6} +
     168\,\I_4^3\,{\I_6}
     + 9\,\I_2^3\,\I_6^2 +
     1242\,{\I_2}\,{\I_4}\,\I_6^2
     }
     \\_{\qquad
     - 1512\,\I_6^3 -
     270\,\I_2^4\,{\I_{10}} -
     11556\,\I_2^2\,{\I_4}\,{\I_{10}} +
     92016\,\I_4^2\,{\I_{10}} +
     37584\,{\I_2}\,{\I_6}\,{\I_{10}} )
     }
     \\
   _{   +
  36\,r^8\,{\I_{10}}\,( \I_2^4\,\I_4^2 -
     17\,\I_2^2\,\I_4^3 + 16\,\I_4^4 -
     6\,\I_2^3\,{\I_4}\,{\I_6} +
     96\,{\I_2}\,\I_4^2\,{\I_6} +
     9\,\I_2^2\,\I_6^2 - 144\,{\I_4}\,\I_6^2 -
     1350\,\I_2^3\,{\I_{10}} }\\
     _{\qquad
     +
     23544\,{\I_2}\,{\I_4}\,{\I_{10}} -
     54432\,{\I_6}\,{\I_{10}} )
     }
     \\
    _{ +
  r^7\,( \I_2^4\,\I_4^4 -
     2\,\I_2^2\,\I_4^5 + \I_4^6 -
     12\,\I_2^3\,\I_4^3\,{\I_6} +
     12\,{\I_2}\,\I_4^4\,{\I_6} +
     54\,\I_2^2\,\I_4^2\,\I_6^2 -
     18\,\I_4^3\,\I_6^2 -
     108\,{\I_2}\,{\I_4}\,\I_6^3 + 81\,\I_6^4
     }
     \\_{\qquad
     +
     30\,\I_2^5\,{\I_4}\,{\I_{10}} +
     156\,\I_2^3\,\I_4^2\,{\I_{10}} +
     1272\,{\I_2}\,\I_4^3\,{\I_{10}} -
     72\,\I_2^4\,{\I_6}\,{\I_{10}} -
     3672\,\I_2^2\,{\I_4}\,{\I_6}\,{\I_{10}} +
     2448\,\I_4^2\,{\I_6}\,{\I_{10}}
     }
     \\_{\qquad
      +
     7236\,{\I_2}\,\I_6^2\,{\I_{10}} -
     1202364\,\I_2^2\,\I_{10}^2 +
     4167936\,{\I_4}\,\I_{10}^2 )
     }
     \\
     _{     -
  4\,r^6\,{\I_{10}}\,\left( \I_2^6 -
     218\,\I_2^4\,{\I_4} -
     512\,\I_2^2\,\I_4^2 - 5832\,\I_4^3 +
     312\,\I_2^3\,{\I_6} +
     18480\,{\I_2}\,{\I_4}\,{\I_6} - 28152\,\I_6^2 +
     2^43^7 67\,{\I_2}\,{\I_{10}} \right)
     }
     \\
   _{  -
  3\,r^5\,( -5\,\I_2^4\,\I_4^3 +
     19\,\I_2^2\,\I_4^4 - 14\,\I_4^5 +
     42\,\I_2^3\,\I_4^2\,{\I_6} -
     96\,{\I_2}\,\I_4^3\,{\I_6} -
     117\,\I_2^2\,{\I_4}\,\I_6^2 +
     126\,\I_4^2\,\I_6^2 +
     108\,{\I_2}\,\I_6^3
     }
     \\_{\qquad
      + 48\,\I_2^5\,{\I_{10}} -
     906\,\I_2^3\,{\I_4}\,{\I_{10}} +
     372\,{\I_2}\,\I_4^2\,{\I_{10}} -
     6120\,\I_2^2\,{\I_6}\,{\I_{10}} +
     85824\,{\I_4}\,{\I_6}\,{\I_{10}} +
     7589376\,\I_{10}^2 )
     }
     \\
   _{
     -
  2\,r^4\,( \I_2^5\,\I_4^2 -
     110\,\I_2^3\,\I_4^3 +
     109\,{\I_2}\,\I_4^4 -
     6\,\I_2^4\,{\I_4}\,{\I_6} +
     810\,\I_2^2\,\I_4^2\,{\I_6} -
     156\,\I_4^3\,{\I_6} + 9\,\I_2^3\,\I_6^2 -
     1917\,{\I_2}\,{\I_4}\,\I_6^2
     }
     \\_{\qquad
     + 1404\,\I_6^3 +
     594\,\I_2^4\,{\I_{10}} +
     24678\,\I_2^2\,{\I_4}\,{\I_{10}} +
     27216\,\I_4^2\,{\I_{10}} -
     140616\,{\I_2}\,{\I_6}\,{\I_{10}} )
     }
     \\
   _{
       - 9\,r^3\,( 4\,\I_2^4\,\I_4^2 -
     116\,\I_2^2\,\I_4^3 + 31\,\I_4^4 -
     24\,\I_2^3\,{\I_4}\,{\I_6} +
     672\,{\I_2}\,\I_4^2\,{\I_6} +
     36\,\I_2^2\,\I_6^2 -
     1008\,{\I_4}\,\I_6^2 - 24\,\I_2^3\,{\I_{10}} }
     \\
     _{\qquad
     +
     36960\,{\I_2}\,{\I_4}\,{\I_{10}} -
     94464\,{\I_6}\,{\I_{10}} )
     }
     \\
      _{  -
  54\,r^2\,\left( 4\,\I_2^3\,\I_4^2 -
     31\,{\I_2}\,\I_4^3 -
     24\,\I_2^2\,{\I_4}\,{\I_6} +
     108\,\I_4^2\,{\I_6} + 36\,{\I_2}\,\I_6^2 -
     504\,\I_2^2\,{\I_{10}} + 9792\,{\I_4}\,{\I_{10}}
     \right)
     }
     \\
   _{ -
  432\,r\,\left( \I_2^2\,\I_4^2 - \I_4^3 -
     6\,{\I_2}\,{\I_4}\,{\I_6} + 9\,\I_6^2 -
     54\,{\I_2}\,{\I_{10}} \right)
-2^83^6\,{\I_{10}}=0}.
\end{array}
\end{equation}

Note that the coefficient of $r^n$ in this expression is an homogeneous weighted polynomial of degree $2n+10$ in the invariants $\I_2, \I_4, \I_6, \I_8$.

\subsection{Jacobian nullwerte and symmetric equations in genus~two}

We now proceed to specialize the results in sections \ref{matrius-algebraiques} and \ref{Nullwerte-Weierstrass}
for hyperelliptic genus two curves. Our starting point now is the normalized period matrix for the Jacobian of a hyperelliptic  curve  defined over a number field, and our goal is the determination of a symmetric model for the curve.
The  procedure described in section \ref{matrius-algebraiques} provides a basis  of algebraic differential forms for the curve, and then proposition \ref{imatge-canonica} gives formulas for the canonical image of the Weierstrass points of the curve with respect to this basis. We obtain the following simple formulas:

\begin{theorem} \label{teorema-gran-g2}
Let $C$ be a genus 2 curve, with field of moduli $K\subset \C$.
Let $Z\in\HH_2$ be a normalized period matrix for $C$.
 Given two odd 2-torsion points $w_1, w_2\in J(C)[2]$.
\begin{itemize}
\item[a)] For every even 2-torsion point $w_0\in J(C)[2]$, the non-singular matrix
$
\Omega_1(w_1,w_2,w_0;Z):=\frac1{2\pi i \q[w_0]} J[w_1,w_2]
$
is a period matrix of a basis of differential forms on $C$ defined over a finite extension $L$ of $K$.
\item[b)] Let $w$ be a third odd 2-torsion point on $J(C)$. The point cut by the hyperplane
$$
\left(
\begin{array}{ll}
\displaystyle\frac{\partial \theta }{\partial z_{1}}(w) &
\displaystyle\frac{\partial \theta }{\partial z_{2}}(w)
\end{array}
\right) \Omega _{1}(w_1,w_2,w_0;Z)^{-1}\left(
\begin{array}{c}
X_{1} \\
X_{2}
\end{array}
\right)
 =0,
$$
in $\P^1(\C)$ has projective coordinates
$\left( [w_1,w] : [w_2,w]|\right)$
and it is independent of $w_0$.
\item[c)] Let $J(C)[2]^{\mbox{odd}}=\{w_1, w_2, w_3, w_4, w_5, w_6\}$. The ratios
\begin{equation}\label{arrels-j}
\ell_{12j}:=\frac{[w_1,w_j]}{[w_2,w_j]}\in \C\cup\{\infty\}, \quad j=1,\dots,6,
\end{equation}
are algebraic over $K$.
\item[d)] The curve $C$ admits the symmetric model
\begin{equation}
\label{equacio-simetrica}
\CC_{12}:\quad Y^2=X(X-\ell_{123})(X-\ell_{124})(X-\ell_{125})(X-\ell_{126})
\end{equation}
over a finite extension of $K$.
\item[e)]
The symmetric discriminant $\D_{12}$ of the symmetric model $\CC_{12}$ is given by:
$$
\D_{12}=\frac{[w_1,w_2]^{16}}{\left([w_1,w_3][w_1,w_4][w_1,w_5][w_1,w_6]\right)^4}.
$$

\end{itemize}

\end{theorem}

\begin{proof}
{\it a}) The discussion in  section \ref{matrius-algebraiques} shows that $\Omega_1(w_1,w_2,w_0;Z)$ is a period matrix for certain basis  $\eta_1,\eta_2$ of $H^0(J(C),\Omega^1_{/\overline{K}})$  defined over a finite extension of $K$. Since the Abel-Jacobi map $\Pi: C \rightarrow J(C)$ is defined over a finite extension of   $K$, the  forms $\w_1=\Pi^\ast\eta_1,\w_2=\Pi^\ast\eta_2$  are  defined over a finite extension $L$ of $K$, and $\Omega_1(w_1,w_2,w_0;Z)$ is a period matrix for them.

\noindent
{\it b}) From the equality
$$
\Omega_1(w_1,w_2,w_0;Z)^{-1}=\lambda
\left(
\begin{array}{ccc}
\ds \frac{\partial \q[w_2]}{\partial z_2} (0) &  \ds -\frac{\partial \q[w_1]}{\partial z_2}(0) \\ \\
\ds -\frac{\partial \q[w_2]}{\partial z_1} (0) &  \ds \frac{\partial \q[w_1]}{\partial z_1}(0) \\
\end{array}
\right),
$$
with $\lambda\in\C^\ast$, a simple matrix calculation shows that
$$
\begin{array}{c}
\left(
\displaystyle\frac{\partial \theta }{\partial z_{1}}(w),
\displaystyle\frac{\partial \theta }{\partial z_{2}}(w)
\right)
\Omega_1(w_1,w_2,w_0;Z)^{-1}
=\lambda\, {h(w)}\left([w,w_2], [w_1,w]
\right),
\end{array}
$$
with $h(w)\neq0$,
and the assertion follows.

{\it c}) The ratios of Jacobian Nullwerte can be reduced to ratios of differences of $x$-coordinates of the Weierstrass points of any algebraic model $Y^2=f(X)$ of $C$, by theorem \ref{Thomae-J}, and hence they are algebraic themselves.

 {\it d}) By \cite[p. 399]{Goro}, there exists two functions $x,y\in\overline{K}(C)$ such that $\w_1=dx/y$, $\w_2=x\,dx/y$, providing a model $y^2=f(x)$ for $C$, defined over $L$. For this model we know by
\ref{hiperpla} that their Weierstrass points have the ratios $\ell_{12j}$ as $x$-coordinates.
Using theorem \ref{Thomae-J} is easy to see that $\prod_{j=3}^6\ell_{12j}=\pm1$, and hence $\CC_{12}$ is a symmetric model for $C$.

 {\it d}) The expression for $\D_{12}$ follows from the equality:
 $$
 \ell_{12i}-\ell_{12j}=\frac{[w_1,w_2][w_i,w_j]}{[w_i,w_2][w_j,w_2]}.
 $$
\end{proof}

After theorem \ref{teorema-gran-g2}, we have a complete and effective solution to the hyperelliptic Torelli problem in genus~2. We have applied this result in three different situations, to present irreducible abelian surfaces with extra multiplications as Jacobians of curves (\cite{GGG}, \cite{GGR}, \cite{BG}).

Theorem \ref{teorema-gran-g2} also has a number of theoretical consequences: we can rephrase the properties
 of the symmetric roots  in terms of the corresponding expressions with the Jacobian Nullwerte, thus
 providing elementary proofs for relations between them.
For instance:

\begin{proposition}
Let $J(C)[2]^{\mbox{odd}}=\{w_1, w_2, w_3, w_4, w_5, w_6\}$. We have
$$
[w_1,w_3] [w_1, w_4] [w_1,w_5] [w_1,w_6] =\pm [w_2,w_3] [w_2, w_4] [w_2,w_5] [w_2,w_6]
$$
\end{proposition}

This result could be proved by means of the Rosenhain formula, but the interpretation of the quotients $[w_1,w_i]/[w_2,w_i]$ as symmetric roots gives it immediately.

\subsection{Thetanullwerte and symmetric equations in genus~two}
\label{thetanullwerte}

We now combine the expression (\ref{arrels-j}) of the symmetric roots in terms of Jacobian Nullwerte with Rosenhain formula (\ref{Rosenhain}), obtaining expressions for the symmetric roots as quotients of Thetanullwerte:
\begin{proposition}
Let $J(C)[2]^{\mbox{odd}}=\{w_i, w_j, w_k, w_a, w_b, w_c\}$. We have:
\begin{tabular}{rl}
$\ell_{ijk}$
&
$=\pm
\prod_{\begin{array}{c}
r=1,\dots,6\\ r\neq i,j,k
\end{array}}\ds
\frac{\q[w_i+w_k-w_r]}{\q[w_j+w_k-w_r]}
$\\
&
$=\ds\pm\frac{\q[w_i+w_k-w_a]\q[w_i+w_k-w_b]\q[w_i+w_k-w_c]}{\q[w_j+w_k-w_a]\q[w_j+w_k-w_b]\q[w_j+w_k-w_c]} $ \\\\
&
$=\ds\pm\frac{\q[w_i+w_k-w_a]\q[w_i+w_k-w_b]\q[w_i+w_k-w_c]}{\q[w_i+w_a-w_b]\q[w_i+w_a-w_c]\q[w_i+w_b-w_c]};  $ \\ \\
\multicolumn{2}{l}
(The sign depends only on $w_i, w_j, w_k$ and can be explicitly determined.)
\\ \\
$\D_{ij}$&
$
=\frac{\left(\q[w_i+w_j+w_k]\q[w_i+w_j+w_a]\q[w_i+w_j+w_b]\q[w_i+w_j+w_c]\right)^{12}}
{\left(\q[w_i+w_k+w_a]\q[w_i+w_k+w_b]\q[w_i+w_k+w_c]\q[w_i+w_a+w_b]\q[w_i+w_a+w_c]\q[w_i+w_b+w_c]\right)^8}. $
\end{tabular}
\end{proposition}

Weber \cite{Weber} and Takase \cite{Takase} have  given expressions for the roots of a Rosenhain model for $C$ in terms of even Thetanullwerte. Takase's formulas are simpler, since they are quotients of Thetanullwerte. One can recover these formulas from the proposition above, just by deriving a Rosenhain model from our symmetric model.

The simple expressions in the proposition (or Takase's formulas) worth attention for practical applications. In case we want to compute a (symmetric) equation for $C$ from its period matrix, it happens that actually only six different Thetanullwerte are involved in the computation  of a set of  symmetric roots of $C$. In particular, the computation of the Igusa invariants of $C$ by means of these formulas requires only six numerical evaluations of the Theta function. This represents a gain of $40\%$ with respect to the methods applied in   (\cite{Weber}, \cite{Weng}, \cite{Wang}). Moreover, the minimality properties of the symmetric models (theorem \ref{minimal}) suggest that their coefficients should be relatively small.
In case we have some extra information about the arithmetic of the curve (for instance, if we know that its Jacobian variety has complex multiplication and we know over which primes its reduction is not irreducible), we will be able to bound the denominators appearing on the symmetric equations, a crucial point to be sure that numerical results have enough precision to be correct.


\section{Genus 3 curves}

We shall now describe an effective solution to the Torelli problem for hyperelliptic genus 3 curves, based on the combination of Jacobian Nullwerte and symmetric roots.

We are given a normalized period matrix $Z\in\HH_3$, corresponding to a genus 3 hyperelliptic curve $C$, and we are asked for an equation of the curve.

We may suppose that the period matrix $Z$ comes from a model $Y^2=(X-\alpha_1)\dots(X-\alpha_8)$ of $C$. The Weierstrass points of this model are then $W_i:=(\alpha_i,0)$. There are twenty-eight odd two torsion points in $J(C)$, and they are given by the degree 2 divisors of the form $W_i+W_j$. We shall write $w_{rs}=\Pi(W_r+W_s)$.
Note that $w_{rs}+w_{st}=w_{rt}$.

The following result is particular case of corollary \ref{J-alfa} for genus~3 hyperelliptic curves:

\begin{lemma} Let $m,n,r,s,t\in \{1,2,\dots,8\}$.  We have:
$$
\mu_{mnrs}=\frac{[w_{mt},w_{ts},w_{sn}][w_{nt},w_{tr},w_{rn}]}{[w_{mt},w_{tr},w_{rn}][w_{nt},w_{ts},w_{sn}]}.
$$
\end{lemma}

The first step in the computation of a symmetric model for $C$ is the proper identification  of the elements of $J(C)[2]^{\rm odd}$ with the divisors $W_i+W_j$. We proceed as follows: we take an arbitrary pair
 $w_1, w_2 \in J(C)[2]^{\rm odd}$ such that $w_3:=w_1+w_2$ is also odd. This assures that these points come from three divisors geometrically well posed in the sense of proposition \ref{sistema fonamental hipereliptic}:
$$
w_1:=w_{23}=\Pi(W_2+W_3), w_2=w_{13}=\Pi(W_1+W_3),  w_3=w_{12}=\Pi(W_1+W_2),
$$
(a formal re-labelling of the $\alpha_i$ may be necessary).
We now look for the remaining five points $w\in J(C)[2]^{\rm odd}$ such that $w_2+w, w_3+w$ are simultaneously odd; they must came from divisors $W_1+W$. We can write them as
$$
w_{14}=\Pi(W_1+W_4),\quad w_{15}=\Pi(W_1+W_5),\quad, \dots, w_{18}=\Pi(W_1+W_8).
$$
Since $w_{jk}=w_{1j}+w_{1k}$, we are already in position to apply the lemma above to compute all the $\mu_{12rs}$:
\begin{equation}
\label{mu12}
\mu_{12rs}=\frac{[w_{1t},w_{ts},w_{s2}][w_{2t},w_{tr},w_{r2}]}{[w_{1t},w_{tr},w_{r2}][w_{2t},w_{ts},w_{s2}]}.
\end{equation}
There are eighteen 2-torsion points involved in this computation, so that, in principle, we will have to compute 54 theta-derivatives (but these calculations  are highly parallelizable).

Finally, we compute a set of symmetric roots for $C$. After lemma \ref{lema-mu-invariants} we have
$$
\ell_{123}=\sqrt[6]{\prod_{k\neq 1,2,3}\mu_{123k}}
$$
(no matter which root we take), and then
$$
\ell_{12k}=\mu_{123k}^{-1}\ell_{123}, \qquad k=4,\dots,8.
$$
We have obtained:

\begin{theorem}
Let us denote by $[ab,cd,ef]$ the Jacobian Nullwerte $[w_{a,b},w_{b,c},w_{c,d}](Z)$. A set of symmetric roots for $C$ is:
$$
\begin{array}{l}
\ell_{123}=\sqrt[6]{
\frac{[  {14},  {48},  {28}]\,
    [  {18},  {48},  {24}]\,
    [  {18},  {58},  {25}]\,
    [  {18},  {68},  {26}]\,
    [  {18},  {78},  {27}]\,
    [  {24},  {34},  {23}]\,
{[  {28},  {38},  {23}]}^4}{[  {14},  {34},  {23}]\,
    {[  {18},  {38},  {23}]}^4\,
    [  {24},  {48},  {28}]\,
    [  {28},  {48},  {24}]\,
    [  {28},  {58},  {25}]\,
    [  {28},  {68},  {26}]\,
    [  {28},  {78},  {27}]}
}
\\\\
\ell_{124}=\sqrt[6]{
\frac{[  {14},  {48},  {28}]\,
    {[  {18},  {38},  {23}]}^2\,
    [  {18},  {58},  {25}]\,
    [  {18},  {68},  {26}]\,
    [  {18},  {78},  {27}]\,
    [  {24},  {34},  {23}]\,
    {[  {28},  {48},  {24}]}^5}{[  {14},  {34},  {23}]\,
    {[  {18},  {48},  {24}]}^5\,
    [  {24},  {48},  {28}]\,
    {[  {28},  {38},  {23}]}^2\,
    [  {28},  {58},  {25}]\,
    [  {28},  {68},  {26}]\,
    [  {28},  {78},  {27}]}}
    \\\\
\ell_{125}=\sqrt[6]{
\frac{[  {14},  {48},  {28}]\,
    {[  {18},  {38},  {23}]}^2\,
    [  {18},  {48},  {24}]\,
    [  {18},  {68},  {26}]\,
    [  {18},  {78},  {27}]\,
    [  {24},  {34},  {23}]\,
{[  {28},  {58},  {25}]}^5}{[  {14},  {34},  {23}]\,
    {[  {18},  {58},  {25}]}^5\,
    [  {24},  {48},  {28}]\,
    {[  {28},  {38},  {23}]}^2\,
    [  {28},  {48},  {24}]\,
    [  {28},  {68},  {26}]\,
    [  {28},  {78},  {27}]}}
\\
\\
\ell_{126}=\sqrt[6]{\frac{[  {14},  {48},  {28}]\,
    {[  {18},  {38},  {23}]}^2\,
    [  {18},  {48},  {24}]\,
    [  {18},  {58},  {25}]\,
    [  {18},  {78},  {27}]\,
    [  {24},  {34},  {23}]\,
    {[  {28},  {68},  {26}]}^5}{[  {14},  {34},  {23}]\,
    {[  {18},  {68},  {26}]}^5\,
    [  {24},  {48},  {28}]\,
    {[  {28},  {38},  {23}]}^2\,
    [  {28},  {48},  {24}]\,
    [  {28},  {58},  {25}]\,
    [  {28},  {78},  {27}]}}
    \\\\
\ell_{127}=\sqrt[6]{
\frac{[  {14},  {48},  {28}]\,
    {[  {18},  {38},  {23}]}^2\,
    [  {18},  {48},  {24}]\,
    [  {18},  {58},  {25}]\,
    [  {18},  {68},  {26}]\,
    [  {24},  {34},  {23}]\,
    {[  {28},  {78},  {27}]}^5}{[  {14},  {34},  {23}]\,
    {[  {18},  {78},  {27}]}^5\,
    [  {24},  {48},  {28}]\,
    {[  {28},  {38},  {23}]}^2\,
    [  {28},  {48},  {24}]\,
    [  {28},  {58},  {25}]\,
    [  {28},  {68},  {26}]}
} \\\\
\ell_{128}=\sqrt[6]{\frac{{[  {14},  {34},  {23}]}^5\,
    [  {18},  {48},  {24}]\,
    [  {18},  {58},  {25}]\,
    [  {18},  {68},  {26}]\,
    [  {18},  {78},  {27}]\,
    {[  {24},  {48},  {28}]}^5\,
{[  {28},  {38},  {23}]}^4}{{[  {14},  {48},  {28}]}^5\,
    {[  {18},  {38},  {23}]}^4\,
    {[  {24},  {34}  {23}]}^5\,
    [  {28},  {48},  {24}]\,
    [  {28},  {58},  {25}]\,
    [  {28},  {68},  {26}]\,
    [  {28},  {78},  {27}]}
}
\end{array}
$$

\end{theorem}

These expressions are not unique, since they depend on the value of $t$ chosen in (\ref{mu12}) to determine $\mu_{123k}$. In any case, the chance to pick two different $t$ for the same $k$ gives a lot of equalities between quotients of Jacobian Nullwerte. We give just one example
:
\begin{proposition} For every hyperelliptic period matrix $Z\in\HH_3$
$$
\frac{[w_{15},w_{45},w_{24}](Z)\,
    [w_{25},w_{35},w_{23}](Z)}{[w_{15},
     w_{35},w_{23}](Z)\,
    [w_{25},w_{45},w_{24}](Z)}
 =
 \frac{[w_{16},w_{46},w_{24}](Z)\,
    [w_{26},w_{36},w_{23}](Z)}{[w_{16},
     w_{36,}w_{23}](Z)\,
    [w_{26},w_{46},w_{24}](Z)}
$$
\end{proposition}

We finally remark that one can also express the symmetric roots $\ell_{12k}$ as quotients of Thetanullwerte, using Frobenius formula (\ref{Frobenius-g3}). For instance:
$$
\ell_{123}=\sqrt[6]{
\frac{\q[1345]^3\q[1346]\q[1367]\q[1368]\q[1456]^3\q[1478]^5}
{\q[1356]^5\q[1378]^3\q[1457]\q[1458]\q[1578]\q[1678]^3},
}
$$
where  the Thetanullwerte $\q[\Pi(W_a+W_b+W_c+W_d)]$ has been written  $\q[abcd]$.
An important remark for computational purposes is that the whole set of formulas for the symmetric roots involve only twelve even Thetanullwerte.


\begin{thebibliography}{ACGH 84}

\bibitem{ACGH84}
 Arbarello, E., Cornalba, M., Griffiths, P.A.,
Harris, J., \textit{Geometry of algebraic curves}, Grundlehren Math.
Wiss., 267, Springer V., New York, 1985; MR0770932 (86h:14019).




\bibitem{BG} Bayer, P., Gu\`{a}rdia, J.; Hyperbolic uniformization of the Fermat curves,
\textit{Ramanjujan J.}, \textbf{12} (2006), pp. 207-223.

\bibitem{Antwerp} B. J. Birch; W. Kuyk. (eds.),  {\it Modular functions of one variable, IV (Proc. Internat. Summer School,
 Univ. Antwerp, Antwerp, 1972)},  Lecture Notes in Math.,  476, Springer
, Berlin, 1975; MR0376533 (51 \#12708).





\bibitem{Biel} Cardona, G., Quer, J., Field of moduli and field of definition for curves of genus 2,
 in {\it Computational aspects of algebraic curves}, 71--83, World Sci. Publ., Hackensack, NJ, 2005; MR2181874 (2006h:14036).


\bibitem{Cremona}
Cremona, J. E., {\it Algorithms for modular elliptic curves} , Cambridge Univ. Press
, Cambridge, 1992; MR1201151 (93m:11053).


\bibitem{Frobenius}  Frobenius, F.G., \"{U}ber die constanten
Factoren der Thetareihen, \textit{J. reine angew. Math.}, \textbf{98} (1885), pp.
241--260.

\bibitem{Goro} Gonz\'{a}lez Jim\'{e}nez, E., Gonz\'{a}lez-Rovira, J.;
Modular curves of genus 2, \textit{Math.  Comp.}, (2003), no.~241, 397--418; MR1933828 (2003i:11078).



\bibitem{GGG} Gonz\'{a}lez Jim\'{e}nez, E., Gonz\'{a}lez-Rovira, J,, Gu\`{a}rdia, J.; Computations on modular {J}acobian surfaces,
in {\it Algorithmic number theory (Sydney, 2002)}, 189--197, Lecture Notes in Comput. Sci., 2369, Springer, Berlin, 2002; MR2041083 (2005c:11074).

\bibitem{GGR} Gonz\'{a}lez-Rovira, J., Gu\`{a}rdia, J., Rotger, V.; Abelian surfaces of $\mbox{GL}_2$ type as Jacobians of curves, \textit{Acta Arith.} {\bf 116} (2005), no.~3, 263--287; MR2114780 (2005m:11107).




\bibitem{Grant}
Grant, D.,
A generalization of Jacobi's derivative formula to dimension two.
{\textit J. Reine Angew. Math.} {\bf 392} (1988), 125--136; MR0965060 (89m:14024).



\bibitem{Guardia} Gu\`{a}rdia, J.,  Jacobian nullwerte and algebraic equations,
\textit {J. Algebra} {\bf 253} (2002), no.~1, 112--132; MR1925010 (2004a:14032).

\bibitem{Guardia-MRL} Gu\`{a}rdia, J., Jacobi thetanullwerte, periods of elliptic curves and minimal equations,
\textit{Mathematical Research Letters}, {\bf 11} (2004), no.~1, 115--123; MR2046204 (2005b:11053).

\bibitem{GuToVe} Gu\`{a}rdia, J.; Torres, E.; Vela, M., Stable models of elliptic curves, ring class fields and complex multiplication, in {\it Algorithmic number theory}, 250--262, Lecture Notes in Comput. Sci., 3076, Springer, Berlin, 2004; MR2137358 (2005m:11103).






\bibitem{Igusa-Jacobi formula}  Igusa, J.I., On Jacobi's derivative
formula and its generalizations, \textit{Amer. J. Math.}, {\bf 102} (1980), no.~2, 409--446; MR0564480 (82e:14053).

\bibitem{Igusa-odd} Igusa, J. I.,
On the nullwerte of Jacobians of odd theta functions, in {\it Symposia Mathematica, Vol. XXIV (Sympos., INDAM, Rome, 1979)}, 83--95, Academic Press
, London, 1981; MR0619242 (83e:14030).

\bibitem{Igusa-Poincare} Igusa, J. I., Problems on abelian
functions at the time of Poincar\'{e} and some at  present, \textit{
Bull. Amer. Math. Soc. (N.S.)}   {\bf 6} (1982), no.~2, 161--174; MR0640943 (83f:14035).




\bibitem{Igusa-multip}  Igusa, J. I.,  Multiplicity one theorem
and problems related to Jacobi's formula, \textit{
 Amer. J. Math.}  {\bf 105} (1983), no.~1, 157--187; MR0692109 (84h:14052).






\bibitem{Lockhart}  Lockhart, P., On the discriminant of a
hyperelliptic curve, \textit{Trans. Am. Math. Soc.} {\bf 342} (1994), no.~2, 729--752; MR1195511 (94f:11054).




\bibitem{Magma} MAGMA, {\tt http://magma.math.usyd.edu.au/magma/ }, \textit{Universtity of Sydney}, 2004.

\bibitem{Mestre}
Mestre, J.F., Construction de courbes de genre $2$ \`a partir de leurs modules , in {\it Effective methods in algebraic geometry (Castiglioncello, 1990)}, 313--334, Progr. Math., 94, Birkh\"auser, Boston, Boston, MA, 1991; MR1106431 (92g:14022).



\bibitem{McKean-Moll} McKean, H., and Moll, V., \textit{Elliptic curves: function theory, geometry, aarithmetic}, Cambridge Univ. Press, Cambridge, 1997; MR1471703 (98g:14032).




\bibitem{Mumford-Tata2}  Mumford, D., \textit{ Tata lectures on theta, II},
Progress in Mathematics, 43, Birkh\"{a}user Boston, Boston, MA, 1984; MR0742776 (86b:14017).






\bibitem{Rosenhain}  Rosenhain, G., M\'{e}moire sur les fonctions de
deux variables et \`{a} quatre p\'{e}riodes qui sont les inverses des
int\'{e}grales ultra-elliptiques de la premi\`{e}re classe, \textit{%
M\'{e}moires des savants \'{e}trangers}, XI (1851), pp. 362-468.


\bibitem{Shimura-CM} Shimura, G., \textit{Abelian varieties with complex multiplication and modular functions}, Princeton Series, 46, Princeton Univ. Press, Princeton, NJ, 1998; MR1492449 (99e:11076).



\bibitem{Silverman} Silverman, J., {\it The arithmetic of elliptic
curves}, G.T.M., 106, Corrected reprint of the 1986 original,
Springer, New York, 1992; MR1329092 (95m:11054).


\bibitem{Takase} Takase, K., A generalization of Rosenhain's
normal form for hyperelliptic curves with an application,\textit{Proc. Japan Acad. Ser. A Math.
Sci.}, {\bf 72} (1996), no.~7, 162--165; MR1420607 (98a:14043).




\bibitem{Thomae}  Thomae, J., Beitrag zur Bestimmung von $\theta
(0,0,...,0)$ durch die Klassenmoduln algebraischer Funktionen, \textit{J.
reine angew. Math.}, 71 (1870), pp. 201--222.

\bibitem{Wamelen} van Wamelen, P.  Examples of genus two CM curves defined over the rationals. \textit{Math. Comp.} {\bf 68} (1999), no.~225, 307--320; MR1609658 (99c:11079).

\bibitem{Wang} Wang, X. D.,  $2$-dimensional simple factors of $J\sb 0(N)$. \textit{ Manuscripta Math.} {\bf 87} (1995), no.~2, 179--197; MR1334940 (96h:11059).



\bibitem{Weil} Weil, A., Sur les p\'{e}riodes des int\'{e}grales ab\'{e}liennes, \textit{Com. on Pure and Applied Math.}, 29 (1976), pp. 813-819.

\bibitem{Weng} Weng, A., A class of hyperelliptic CM-curves of genus three,
\textit  {J. Ramanujan Math. Soc.}  {\bf 16} (2001), no.~4, 339--372; MR1877806 (2002k:11099).

\bibitem{Weber} Weber, H.-J., Hyperelliptic simple factors of $J\sb 0(N)$ with dimension at least $3$.
\textit{Experiment. Math.}   {\bf 6} (1997), no.~4, 273--287; MR1606908 (99e:14054).







\end{thebibliography}
\end{document}